\documentclass[12pt]{article}
\usepackage[latin1]{inputenc}
\usepackage{amssymb}
\usepackage{amsmath}
\usepackage{amsthm}
\usepackage{tikz}

\newcommand{\C}{\mathcal{C}}
\newcommand{\F}{\mathbb{F}}
\newcommand{\R}{\mathbb{R}}

\newcommand{\Z}{\mathbb Z}

\newcommand{\Fix}{\mathrm{Fix}}

\theoremstyle{definition}

\newtheorem{thm}{Theorem}[section]

\newtheorem{prop}[thm]{Proposition}

\newtheorem{lem}[thm]{Lemma}
 
\newtheorem{rem}[thm]{Remark}

\newtheorem{defi}[thm]{Definition}

\newtheorem{cor}[thm]{Corollary}

\newtheorem{claim}[thm]{Claim}

\newcommand{\vs}{\vspace{0.3cm}}

\newcommand{\vsp}{\vspace{0.2cm}}

\date{}
\author{}

\begin{document}

\title{Orderings and flexibility of some subgroups of $Homeo_+(\R)$}
\author{Juan Alonso, Joaqu\'\i n Brum\footnote{The second author was supported by ANII PhD Scolarship.} and  Crist\'obal Rivas\footnote{The third author acknowledges the support of CONICYT via FONDECYT 1150691 and via PIA 79130017.}}
\maketitle
\abstract{In this work we exhibit flexibility phenomena for some (countable) groups acting by order preserving homeomorphisms of the line. More precisely, we show that if a left orderable group  admits an amalgam decomposition of the form $G=\F_n*_{\Z} \F_m$ where $n+m\geq 3$, then every faithful action of $G$ on the line by order preserving homeomorphisms can be approximated by another action (without global fixed points) that is {\em not} semi-conjugated to the initial action. We deduce that $\mathcal{LO}(G)$, the space of left orders of $G$, is a Cantor set.

In the special case where $G=\pi_1(\Sigma)$ is the fundamental group of a closed hyperbolic surface, we found finer techniques of perturbation. For instance, we exhibit a single representation whose conjugacy class in dense in the space of representations. This entails that the space of representations without global fixed points of $\pi_1(\Sigma)$ into $Homeo_+(\R)$ is connected, and also that the natural conjugation action of $\pi_1(\Sigma)$ on $\mathcal{LO}(\pi_1(\Sigma))$ has a dense orbit.

}

\vsp\vsp

\noindent{\bf MSC 2010 classification}: 20F16, 22F05, 37C85, 37F15.

\section{Introduction} \label{sec intro}

Given a group $G$, the space  $Rep(G, Homeo_+(M))$ of representations of $G$ into the group of orientation preserving homeomorphisms of an orientable manifold $M$ is a classical object that encompasses many different areas of mathematics (see for instance \cite{GRGA}). When the manifold $M$ has dimension one (that is $M$ is either $\R$ or $S^1$), 
allowing a faithful action on $M$ has an algebraic counterpart in terms of left invariant (linear or circular) orders  \cite{Calegari, clay rolfsen, GOD}. For instance, from any faithful action on the line of a group $G$ one can {induce} a total left invariant linear order on $G$ ({\em left order} for short), and conversely, from any left order on a (countable) group $G$ one can produce a faithful action of $G$ on the line by orientation preserving homeomorphisms which is unique up to conjugation. This is the so called {\em dynamical realization} of the left order. See \S \ref{sec orderings} for details.

The counterpart of $Rep(G, Homeo_+(\R))$ is the  {\em space of left orders on $G$}, here denoted $\mathcal{LO}(G)$, which is the set of all left orders of $G$ endowed with a natural topology that makes it a Hausdorff, totally disconnected and compact space \cite{sikora}. Linnell showed that this space is either finite or uncountable \cite{linnell}, hence contains a Cantor set when it is infinite. Although of very different nature (one is {\em continuous} while the other is totally disconnected) there are some relationships between these two spaces. For instance it is implicit in Navas's \cite{navas orders} that the dynamical realization of an isolated left order of $G$ is a   {\em locally rigid} action of $G$ on the line, meaning that any sufficiently close representation is {\em semi-conjugated} to it (see \S\ref{sec flexibility rigidity} for definitions, and Proposition \ref{prop rigidity} for an explicit proof). In fact, recently in \cite{mann rivas} a complete characterization of isolated orders in terms of a strong form of rigidity was obtained. In this work, we focus on a closely related question:
$$\text{\em Is there an algebraic characterization of groups allowing isolated left orders?}$$

To this day we only count with some partial results. Tararin obtained an algebraic classification of groups allowing only finitely many left orders, see \cite[Theorem 5.2.1]{KM}. They all turn out to be solvable, and in \cite{rivas tessera} a classification of (virtually) solvable groups allowing isolated left orders was obtained: they all fit in Tararin's classification. On the other hand, free groups \cite{mccleary, navas orders} and more generally free products of left orderable groups \cite{rivas_free} admit no isolated left orders, whereas -for instance- $\F_2\times \Z$ has infinitely many conjugacy classes of isolated left orders \cite{mann rivas}. In this work we generalize the result for free products to groups allowing certain decompositions as amalgamated free products.

Let $G=\F_n*_{w_1=w_2}\F_m$, with $n\geq 2$, $m\geq 1$, be the amalgamated product of the free groups $\F_n$ and $\F_m$ identifying the cyclic subgroups $\langle w_1\rangle\subseteq\F_n$ and $\langle w_2\rangle\subseteq\F_m$.  We will deduce our result about orders of $G$ from a result of {\em flexibility} of its representations. Roughly speaking, a representation is {\em flexible} if it can be approximated by fixed-point-free representations that are not semi-conjugated to it (see Definition \ref{def flexibility}).

\begin{thm} \label{teo amalgama} {\em Let $G=\F_n*_{w_1=w_2}\F_m$ with $n\geq 2$, $m\geq 1$. Then, every representation  $\rho:G\to Homeo_+(\R)$ without global fixed points is flexible.}
\end{thm}

Proposition \ref{prop rigidity} then entails

\begin{cor}  {\em Let $G=\F_n*_{w_1=w_2}\F_m$ with $n\geq 2$. The space of left orders of $G$ has no isolated points.}

\end{cor}

Some comments are in order. The first one is that a group $G$ as in Theorem \ref{teo amalgama} above is always left orderable since for this it suffices to have an order preserving isomorphism between $\langle w_1\rangle$ and $\langle w_2\rangle$, see \cite{BG camb} (alternatively, they are one-relator and torsion free, hence left orderable \cite{brodskii}). For a general condition for orderability of amalgams see \cite{BG lond}.  Secondly, prior to this work the amalgamated free product has been used to construct groups having and infinite space of left orders that {\em contains} isolated points, such as the groups $\langle a,b\mid a^n=b^m\rangle$ from \cite{ito alg, navas hecke} and the groups constructed by Ito's iterative methods \cite{ito tohoku, ito ggd}. In particular the condition $n+m\geq 3$ in Theorem \ref{teo amalgama} is sharp. Last but not least, Theorem \ref{teo amalgama} should be compared with the work of Mann \cite{mann}, where she shows that $G=\pi_1(\Sigma)$ the fundamental groups of an orientable, closed, hyperbolic surface, has special representations into $Homeo_+(S^1)$ (the so called {\em geometric representations}) which are {\em fully rigid}, meaning that their connected component inside $Rep(\pi_1(\Sigma), Homeo_+(S^1))$ is made of a single semi-conjugacy class. By contrast, Theorem \ref{teo amalgama} implies that the semi-conjugacy class of any representation of $\pi_1(\Sigma)$ (and more generally any group as in Theorem \ref{teo amalgama}) into $Homeo_+(\R))$ has empty interior.


Theorem \ref{teo amalgama} is deduced from a technical lemma involving perturbations of representations of the free group into $Homeo_+(\R)$ under some conditions on the image of a specific element $w\in \F_n$ (see Lemma \ref{lem amalgama}). In the special case where  $w$ is a commutator of some generators of $\F_n$, we obtain finer perturbations techniques (see Lemma  \ref{lem commutator semicont} and Lemma \ref{lem commutator}) which allow us to show that $\mathcal{LO}(G)$ is a Cantor set whenever $G$ is a countable left orderable group allowing a decomposition of the form $H*_{h=[a,b]}\F_2(a,b)$. See Theorem \ref{teo manija}. But they provide much more! Indeed, when we restrict our attention to $G=\pi_1(\Sigma)$ we found other results with a strong flexible flavor. For instance, we show how to perturb a given representation in order to blow up global fixed points. Precisely we prove

\begin{thm}\label{teo bonatti} {\em Let $M$ be $\R$ or $S^1$, $\rho\in Rep(\pi_1(\Sigma), Homeo_+(M))$ and $U$  a neighbourhood of $\rho$. Then, there exists $\rho'\in U$ having no global fixed points. }

\end{thm}

We also show how to build a very special representation of $\pi_1(\Sigma)$ into $Homeo_+(\R)$, whose existence can be thought of as a strong form of flexibility. Actually, results such as Theorem \ref{teo conjugacy} and Theorem \ref{teo dense orbit} below, were only known to hold for the non-Abelian free groups \cite{clay, rivas_free}.

\begin{thm} \label{teo conjugacy}  {\em There is a representation of $\pi_1(\Sigma)$ into $Homeo_+(\R)$ without global fixed points, whose conjugacy class under $Homeo_+(\R)$ is dense in $Rep(\pi_1(\Sigma),Homeo_+(\R))$.} 

\end{thm}

Since conjugacy classes are path connected we immediately  obtain (compare with Remark \ref{rem alexander})
\begin{cor} {\em The space of representation without fixed points $Rep_\#(\pi_1(\Sigma),Homeo_+(\R))$} is connected.
\end{cor}

The counterpart of Theorem \ref{teo conjugacy} in the context of group orders is the theorem below. Recall that a group acts on its space of left orders by conjugation, see \S \ref{sec orderings}.

\begin{thm} \label{teo dense orbit} {\em There is a left order on $\pi_1(\Sigma)$ whose orbit under the natural conjugation action is dense in $\mathcal{LO}(\pi_1(\Sigma))$.}

\end{thm}

In \S \ref{sec reduction} we  state our three main lemmas and deduce all the theorems stated in this introduction from them. Lemma \ref{lem amalgama} is proved in \S \ref{sec proof lem amalgama} whereas Lemma \ref{lem commutator semicont} and Lemma \ref{lem commutator} are proved in \S \ref{sec proof commutator}. The preliminary knowledge and definitions to carry out our study is given in \S \ref{sec flexibility rigidity}.

\begin{rem} After our first draft was released it was pointed out to us that Bonatti and Firmo had proved Theorem \ref{teo bonatti} in the category of $C^{\infty}$ diffeomorphisms (Théorème 5.4 in \cite{BF}). Though their techniques are very similar to ours, we choose to provide self-contained proof of Theorem \ref{teo bonatti} for the sake of completeness.
\end{rem}


\vsp

\noindent \textbf{Acknowledgements.} 
The authors thank Kathryn Mann and Andrés Navas for their feedback and interest in this work.

\section{Preliminaries}
\label{sec flexibility rigidity}
\subsection{Flexibility and local rigidity in $Homeo_+(\R)$}

Throughout this work we will deal with the notion of local rigidity. To state it we first need to recall the definition of semi-conjugacy. We say that a non-decreasing map $c:\R \to \R$ is {\em proper} if $c^{-1}$ maps compact sets into bounded sets\footnote{Please note that our definition of proper map is not the traditional one demanding that inverse image of compact sets are compact. If fact, a paradigmatic example  that we want to consider as {\em proper} is the map $c:x\mapsto \max\{n\in\Z\mid n\leq x\}$. For this map we have that $c^{-1}(0)=[0,1)$.}. Note that this is equivalent to demand that the non-decreasing map $c$ satisfies that $c(\R)$ is unbounded in both directions of the line.

\begin{defi} We say that two representations $\rho_i: G\to Homeo_+(\R)$, $i=1,2$, are semi-conjugated if there is a {\em monotone} (i.e. non-decreasing) map $c:\R\to \R$ which is {\em proper} and satisfies \begin{equation}\label{eq semiconjdef}c\circ \rho_1(g)=\rho_2(g)\circ c \qquad \mbox{for all } g\in G.\end{equation}

\end{defi}

Traditionally ({\em e.g.} in \cite{navas-book}), one also insists on the continuity of $c$ above. This has been a pity since that condition causes more inconveniences than the ones it solves. For instance without the continuity assumption one has

\begin{prop} {\em Semi-conjugacy is an equivalence relation.}

\end{prop}

\noindent{\bf Proof:} Reflexivity is obvious, transitivity is natural, we check symmetry.

Suppose (\ref{eq semiconjdef}) holds. Since $c$ is proper we can make $c'(x):=\sup c^{-1}(x)$. Since $c$ is  monotone we have $$c'(x)=\sup c^{-1}((-\infty,x])=\sup\{y\mid c(y)\leq x\}.$$
Since the above supremum is taken over a larger set, monotonicity of $c'$ follows. To show that $c'$ is proper, observe that $c\circ c'(x)=x$ for every point $x$ on which $c$ is continuous. Further, since $c$ is monotone, we have that $c$ is continuous except maybe on a countable subset of the line. It follows that $c\circ c'(\R)$ is unbounded in both directions of the line. In particular, $c'$ is proper.

Finally, covariance also follows since we have
\begin{eqnarray*}
\rho_1(g)(c'(x))
&=& \sup \{\rho_1(g) (y)\mid c(y) \leq x\}\\
&=& \sup\{ z\mid c(\rho_1(g)^{-1}(z))\leq x\} \\
&=& \sup\{ z\mid \rho_2(g)^{-1}(c(z))\leq x \} \\
&=& \sup\{ z\mid c(z)\leq \rho_2(g)(x)\}\\
&=& c'(\rho_2(g) (x)).
\end{eqnarray*}
$\hfill\square$

\vsp

We now let $G$ be a countable and discrete group and $M$ a locally compact oriented manifold  (for our purpose it is enough to consider $M$ as being the real line or the circle). The set $Rep(G,Homeo_+(M))$, of group representations from $G$ to $Homeo_+(M)$, is endowed with the pointwise convergence. That is, $\rho_n$ converges to $\rho$ if and only if $\rho_n(g)$ converges to $\rho(g)$ for all $g\in G$, where the convergence $\rho_n(g)\to\rho(g)$ is given by the compact open topology: for every $\varepsilon$ and for every compact set $K\subset M$ there is $n_0$ such that $n\geq n_0$ implies  $$\sup_{x\in K} |\rho_n(g)(x)-\rho(g)(x)|\leq \varepsilon.$$


Given $\rho\in Rep(G,Homeo_+(M))$ we define $$Fix(\rho)=\{x\in M:\rho(g)(x)=x, \ \forall g\in G\},$$ the set of global fixed points of $\rho$. The subset of representations without global fixed points, with the inherited topology, will be denoted by $Rep_\#(G,Homeo_+(M))$. In this work we will be mainly interested in understanding rigidity inside the space $Rep_\#(G,Homeo_+(\R))$.

\begin{defi} \label{def flexibility}
We say that $\rho\in Rep_\#(G,Homeo_+(\R))$ is {\em locally rigid}, if there is a neighbourhood $U$ of $\rho$ such that every $\rho'\in U \cap Rep_\#(G,Homeo_+(\R))$ is semi-conjugated to $\rho$. If $\rho$ is not locally rigid, then we say that $\rho$ is {\em flexible}.
\end{defi}

\begin{rem} \label{separability} Observe that the convergence of $\rho_n\to \rho$ in $Rep(G, Homeo_+(\R))$ is equivalent to require that $\rho_n(g)\to\rho(g)$ for every $g$ in a generating set of $G$. In particular, if we have a finite generating set for $G$ with $k$ elements, then $Rep(G,Homeo_+(\R))$ is homeomophic to a subset of $Homeo_+(\R)^k$. Since $Homeo_+(\R)$ is metrizable and separable (see for instance \cite{McCoy}), it satisfies the second axiom of countability, and so does $Homeo_+(\R)^k$ and any of its subsets (such as   $Rep(G,Homeo_+(\R))$ and $Rep_\#(G,Homeo_+(\R))$).
\end{rem}

\begin{rem}\label{rem alexander}
Inside the space of all representation, $Rep(G,Homeo_+(\R))$, the so called {\em Alexander trick} can be performed both to retract the space to the trivial representation and/or to find non-semi-conjugated representations arbitrarily close to a given one. Indeed, for  a representation $\rho:G\to Homeo_+(\R)$, we can consider $f_t:\R\to \R$ a continuous path of homeomorphisms over its images with $f_0(x)$ the identity map, and $f_1$ a constant map, and construct
$$\rho_t(g)(x)=\left\{\begin{array}{cl}  f_t\rho(g)f_t^{-1}(x) & \text{if $x\in f_t(\R)$ }\\ x & \text{otherwise.}\end{array}\right. $$
This kind of tricks are not possible inside $ Rep_\#(G,Homeo_+(\R))$. For instance, the space $Rep_\#(\Z,Homeo_+(\R))$ is not connected since the subset of representations satisfying $\rho(a)(x)>x$ for all $x$ (where $a$ is the generator of $\Z$) is open and closed in that space. A similar argument applies for groups of the form $G=\langle a,b|a^m=b^n\rangle$. 
\end{rem}

\subsection{Conjugacy classes in $Homeo_{+}(\R)$.}
\label{sec conjugacy classes}

An important ingredient for proving our results involving commutators, is the description of conjugacy classes in $Homeo_+(\R)$. Luckily, detecting when two given homeomorphism of the line are conjugated is an easy task: {\em it is all encoded in the combinatorics of the  homeomorphisms.} More precisely, if $\psi\phi_1\psi^{-1}=\phi_2$, then $\psi$ maps bijectively the sets
\begin{equation}
  \label{eq comb inf}
  \begin{aligned}
    Fix(\phi_1)=\{x\mid \phi_1(x)=x\}\ &\longleftrightarrow Fix(\phi_2) ,\\
    Inc(\phi_1)=\{x\mid \phi_1(x)>x\} &\longleftrightarrow Inc(\phi_2),\\
    Decr(\phi_1)=\{x\mid \phi_1(x)<x\} &\longleftrightarrow Decr(\phi_2),
  \end{aligned}
\end{equation}
respectively. In fact, $(\ref{eq comb inf})$ characterizes when two homeomorphisms $\phi_1$ and $\phi_2$ are conjugated. 
If there is $\psi\in Homeo_+(\R)$ 
which maps bijectively $Fix(\phi_i)$, $Inc(\phi_i)$ and $Decr(\phi_i)$ ($i=1,2$),
then there exist $\bar{\psi}$ such that  $\bar{\psi}\phi_1\bar{\psi}^{-1}=\phi_2$. This motivates our next

\begin{defi} For $\psi,\phi_1, \phi_2$ homeomorphisms of the real line, we will say that $\psi$ is a {\em weak-conjugation} from $\phi_1$ to $\phi_2$ if \begin{itemize}
\item  $\psi(Fix(\phi_1))=Fix(\phi_2)$ and
\item $\psi(Inc(\phi_1))=Inc(\phi_2)$. \end{itemize}
Additionally, if for an interval $I$ we have that $\psi\phi_1(x)=\phi_2\psi(x)$ for all $x\in I$ we will say that the weak conjugation $\psi$ is {\em strong} on $I$.
\end{defi}

Observe that conjugacy and weak-conjugacy classes are identical, but it is much easier to find/build weak conjugations rather than true conjugating elements. In order to pass from a weak conjugation to a conjugation, the following lemma (and its proof) will be useful. In its proof and throughout the text,  the restriction of a function $\phi$ to a set $C$ will be denoted by $\phi_{|C}$.


\begin{lem}\label{promotion} {\em Let $\psi,\phi_1,\phi_2\in Homeo_+(\R)$. If $\psi$ is a weak-conjugation from $\phi_1$ to $\phi_2$ that is strong on a interval $I$, then there exists a conjugation $\bar{\psi}$ from $\phi_1$ to $\phi_2$ such that:
\begin{itemize}
\item $\bar{\psi}(x)=\psi(x)$ for every $x\in I$ and
\item $\bar{\psi}(x)=\psi(x)$ for every $x\in Fix(\phi_1)$.
\end{itemize}
Moreover, $\bar \psi$ agrees with $\psi$ over $I\cup \phi_1(I)$.}

\end{lem}

\vsp

\noindent {\bf Proof:} We will prove the lemma for the case in which $I=[u,v]$ is compact. The non-compact case is similar.

Since $\psi$ is a weak conjugacy, every connected component $C$ of $\R-(Fix(\phi_1))$ is sent by $\psi$ to a connected component $D$ of $\R-(Fix(\phi_2))$, and  $C\subseteq Inc(\phi_1)$ if and only if $D\subseteq Inc(\phi_2)$. We will define a conjugation on a component $C$ of $\R-(Fix(\phi_1))$. Choose a point $p\in C$. We assume that $C\subseteq Inc(\phi_1)$, as the other case is analogous. Let $J=[p,\phi_1(p))$, $K = [\psi(p),\phi_2(\psi(p)))$, and take $\alpha:J\to K$ an orientation preserving homeomorphism. Notice that $C=\bigcup _{n\in\Z} \phi_1^n(J)$ and $D:=\psi(C)=\bigcup _{n\in\Z} \phi_2^n(K)$.

Define, for $x\in C$, $\psi_C(x)=\phi_2^{-m}(\alpha(\phi_1^m(x))$, where $m$ is the only integer such that $\phi_1^m(x)\in J$. Then $\psi_C$ is a homeomorphism between $C$ and $D$ that conjugates ${\phi_1}_{|C}$ and ${\phi_2}_{|D}$. Defining $\psi_0\equiv \psi$ on $Fix(\phi_1)$, and $\psi_0\equiv \psi_C$ on each component $C$ of  $\R-\Fix(\phi_1)$ (for some choice of $p$ and $\alpha$) gives a conjugation from $\phi_1$ to $\phi_2$. However, $\psi_0$ may not agree with $\psi$ over $I$. To solve this problem, on each component $C$ that intersects $I$, we choose $p\in C$ such that the corresponding $J=[p,\phi_1(p))$ intersects $I$ maximally (that is, $J\cap I$ is either $J$ or $I$), and we choose $\alpha$ so that it agrees with $\psi$ over $J\cap I$. Since $\psi$ is strong on $I$, this $\psi_0$ is a conjugation from $\phi_1$ to $\phi_2$ that agrees with $\psi$ over $I$.

To show the final claim, take $y=\phi_1(x)$ with $x\in I$. Then, $\psi(y)=\psi(\phi_1(x))=\phi_2(\psi(x))=\phi_2(\bar \psi(x))=\bar\psi(\phi_1(x))=\bar\psi(y)$
$\hfill\square$

\subsection{Group orders and dynamical realizations}
\label{sec orderings}

Recall that  a left order on a group $G$ is a total order $\preceq$ satisfying  that given $f,g,h\in G$ such that $f\preceq h $ then $ gf \preceq gh$. If $G$ admits a left order, then we say that $G$ is {\em left orderable}. The reader unfamiliar with this notion may wish to consult \cite{clay rolfsen, GOD, KM}.

A natural topology can be defined on the set of all left orders on $G$, here denoted $\mathcal{LO}(G)$, making it a compact and totally disconnected space. In this topology, a local base at a left order $\preceq\in \mathcal{LO}(G)$ is given by the sets $$V_{g_1,\ldots,g_n}:=\{\preceq' \in\mathcal{LO}(G)\mid id \prec'g_i\},$$ 
where $\{g_1,\ldots,g_n\}$ runs over over all finite subsets of $\preceq$-positive elements of $G$. In particular, a left order is isolated in $\mathcal{LO}(G)$ if there is a finite set $S\subset G$ such that $\preceq$ is the only left order satisfying $$id \preceq s\, , \text{ for every $s\in S$}.$$
When the group is countable this topology is metrizable \cite{clay rolfsen, GOD,sikora}. For instance, if $G$ is finitely generated, and $B_n$ denotes the ball of radius $n$ with respect to a finite generating set, then we can declare that $dist(\preceq_1,\preceq_2)=1/n$, if $B_n$ is the largest ball in which $\preceq_1$ and $\preceq_2$ coincide.

There is also a natural action of a group $G$ on the space $\mathcal{LO}(G)$ by {\em conjugation} of the orders. Precisely, if $\prec$ is a left order on $G$ and $g\in G$, we can define the order $\prec_g$ by $$ h\prec_g k \ \Leftrightarrow \ ghg^{-1}\prec gkg^{-1}.$$
This $\prec_g$ is the result of acting on $\prec$ by $g$, and it is easy to check that this defines a left action by homeomorphisms of $\mathcal{LO}(G)$.

When the group $G$ is countable, for every left order $\preceq$ on $G$, one can attach a fixed-point-free action $\rho: G\to Homeo_+(\R)$ that models the left translation action of $G$ on $(G,\preceq)$, in the sense that
\begin{equation} f\prec g \Leftrightarrow \rho(f)(0)< \rho(g)(0). \label{eq dyn real}\end{equation}
This is the so called, {\em dynamical realization} of $\preceq$ (which is unique up to conjugation), and $0$ is sometimes called the {\em base point}, see \cite{clay rolfsen, GOD, ghys}.

The action of $G$ by conjugation can also be expressed nicely in terms of dynamical realizations. If $\rho$ is a dynamical realization of $\prec$, then $$h\prec_g k \ \Leftrightarrow \ \rho(h)\rho(g)^{-1}(0)<\rho(k)\rho(g)^{-1}(0) .$$
So, a dynamical realization of $\prec_g$ is the conjugation of $\rho$ by $\rho(g)$. Alternatively, one can  see the order $\prec_g$ as the order induced by the representation $\rho$, but ``based" at the point $\rho(g)^{-1}(0)$. 

So far, only two techniques are know to approximate a given left order $\preceq$ on a group $G$. One is to approach it by its own conjugates $(\preceq_g)_{g\in G}$, see for instance \cite{navas orders, rivas tessera}, and the other one, implicit in \cite{navas orders, rivas_free}, is by showing that the dynamical realization of $\preceq$ is not locally rigid. Indeed we have

\begin{prop}\label{prop rigidity} {\em Let $G$ be a left orderable group and $\prec\in \mathcal{LO}(G)$ an isolated order. Then its dynamical realization $\rho\in Rep_{\#}(G,Homeo_{+}(\R))$ is locally rigid. }
\end{prop}

\noindent {\bf Proof:} Take $F\subseteq G$ a finite set so that $\prec$ is the only left order on $G$ satisfying  $id \preceq f$ for all $f\in F$. Let $\rho$ be a dynamical realization of $\prec$. Then, there is a neighbourhood $U\subset Rep_{\#}(G,Homeo_{+}(\R))$ of $\rho$ so that for $\rho'\in U$ and for every $f\in F\setminus\{id\}$ we have $0<\rho'(f)(0)$ if and only if $0<\rho(f)(0)$. Let $\prec'$ be the partial left order defined by $$g_1\preceq' g_2 \text{ if and only if } \rho'(g_1)(0)\leq \rho'(g_2)(0).$$

Since $Stab_{\rho'(G)}(0)$ is left orderable, we can extend the partial order $\preceq'$ to a total left order, that we still call $\prec'$. See for instance \cite[\S2.1]{GOD}. As $\prec'$ agrees with $\prec$ on $F$, we must have $\prec'=\prec$. In particular, this means that $Stab_{\rho'(G)}(0)$ is trivial, since every non-trivial left orderable group has at least two different orders. 

Therefore we have that $\rho'(g_1)(0)<\rho'(g_2)(0)$ if and only if $\rho(g_1)(0)<\rho(g_2)(0)$ for every $g_1, g_2\in G$. Let $\mathcal{O}$ and $\mathcal{O}'$ be the orbits of $0$ under $\rho$ and $\rho'$ respectively. Then $\rho(g)(0)\mapsto \rho'(g)(0)$ is a monotone and $G$-equivariant map, that we call  $\varphi:\mathcal{O}\to \mathcal{O}'$. It can be extended to a semi-conjugacy $c:\R\to\R$ between $\rho$ and $\rho'$ by setting $$c(x) = \sup\{\varphi(y): y\in \mathcal{O}, y\leq x\}. $$
Indeed, the monotone map $c$ is proper because both representations have no global fixed points. The covariance also follows since
\begin{eqnarray*}
\rho'(g)(c(x))
&=& \rho'(g)(\sup\{\varphi(y):y\in\mathcal{O},y\leq x\}) \\
&=& \sup\{\rho'(g)(\varphi(y)):y\in\mathcal{O}, y\leq x\} \\
&=& \sup\{\varphi(\rho(g)(y)): y\in\mathcal{O},y\leq x\} \\
&=& \sup\{\varphi(z):z\in\mathcal{O},z\leq\rho(g)(x)\}\\
&=& c(\rho(g)(x)).
\end{eqnarray*}
$\hfill\Box$

We refer the reader to \cite{mann rivas} for more about orders and rigidity.



\section{Reduction to main lemmas}
\label{sec reduction}

Let's fix some notation. If $\Gamma$ is a group and $S\subseteq \Gamma$, we shall denote by $\langle S\rangle$ the subgroup generated by $S$. The set of elements of $\langle S\rangle$ that can be expressed as multiplications of at most $n$ elements in $S\cup S^{-1}$ is denoted by $B_n(\langle S\rangle)$. We will work with $\Gamma=Homeo_+(\R)$. In this case let $Fix(\langle S\rangle)$ be the set of global fixed points of the subgroup, which is the same as the common fixed points of the elements of $S$. If $G$ is a group and $\rho\in Rep(G,Homeo_+(\R))$ then $Fix(\rho)=Fix(\rho(G))$.

We will deduce the theorems announced in \S \ref{sec intro} from three technical lemmas involving the level sets of the {\em word map}. The proof of these lemmas is postponed to \S \ref{sec proof lem amalgama} and \S \ref{sec proof commutator}.

 The word map associated to $w\in\F_n$ sends each  $\rho\in Rep(\F_n,Homeo_{+}(\R))$ to $\rho(w)\in Homeo_{+}(\R)$. Our first lemma can be seen as a weak form of semi-continuity of the level sets of the word map for a general $w\in\F_n$  under the mild dynamical assumption that $\rho(w)$ does not fix a neighborhood of $\infty$. Precisely, in \S \ref{sec proof lem amalgama} we show

\begin{lem} \label{lem amalgama} {\em Let $w\in\F_n$ be a cyclically reduced word, $\rho\in Rep(\F_n,Homeo_{+}(\R))$ ($n\geq 2$) and $p\in\R$. Then, there exists $q_0\in\R$ such that for every $q>q_0$ with $\rho(w)(q)\neq q$, and every $h\in Homeo_{+}(\R)$ that satisfies \begin{itemize}
\item $h_{|(-\infty,q)}=\rho(w)_{|(-\infty,q)}$

\item $Fix(h)\cap(q,+\infty)$ consists of at most one point, and in case  $h$ fixes a point in $(q,\infty)$, then (the graph of) $h$ transverses the diagonal at that point,
\end{itemize}
there exists $\rho^*\in Rep(\F_n,Homeo_{+}(\R))$ such that \begin{itemize}
\item $\rho^*(x_i)_{|(-\infty,p)}=\rho(x_i)_{|(-\infty,p)}$
\item $\rho^*(w)=h$
\item $Fix(\rho^*)\subseteq Fix(\rho)$.
\end{itemize}}
\end{lem}

The next two lemmas, that will be proved in \S \ref{sec proof commutator}, are for $w=[a,b]=aba^{-1}b^{-1}\in\F_2$. In  this case, the level sets of the word map correspond to  $$ \mathcal{V}_h = \{(f,g)\in Homeo_{+}(\R): [f,g]=h\}.$$
In this case we obtain a stronger version of Lemma \ref{lem amalgama}, as we do not require any (dynamical) condition on $[f,g]$.
\begin{lem} \label{lem commutator semicont} {\em Let $f$, $g$ be homeomorphisms, let $h=[f,g]$ and let $K$ be a proper closed  interval. Then, for all $h'$ coinciding with $h$ over the convex closure of $K\cup f(K)$, there is $(f',g')\in \mathcal V_{h'}$ such that \begin{enumerate}
\item $(f',g')$ coincides with $(f,g)$ over $K$
\item $Fix(\langle f',g'\rangle)$ is contained in $Fix(\langle f,g \rangle)$. Moreover, $Fix(\langle f',g'\rangle)$ is contained in $K\cup f(K)$.
\end{enumerate}}
\end{lem}

Our final lemma says it is possible to perturb a representation inside a fixed $\mathcal{V}_h$, changing its semi-conjugacy class, provided it satisfies the following condition: 

If $K$ is a proper closed interval of the line, we say that a pair of homeomorphisms $(f,g)$ satisfies condition $(*_K)$ if the following holds:
$$ \text{\em $(*_K)$}\hspace{1.8cm} \begin{array}{c}\text{\em There is a point $p_{(f,g,K)}$, not fixed by $[f,g]$, }\\
\text{\em that is outside $u(K)$ for every $u\in B_2(\langle f,g \rangle)$.} \end{array}\hspace{2cm}$$

\begin{lem} {\em \label{lem commutator}Let $K=(-\infty , k]$ be a closed proper interval of $\R$, and let $(f,g)$ be a pair satisfying $(*_K)$. Let $h=[f,g]$. Then  we can choose  $(f_1,g_1)$ and $(f_2,g_2)$ in $\mathcal V_h$ agreeing with $(f,g)$ over $K$, such that \begin{enumerate}
	\item $Fix(\langle f_i,g_i\rangle)$ is contained in $Fix(\langle f,g\rangle)$, for $i=1,2$
	\item  $g_1(x)> x$ and $g_2(x)< x$ for $x$ large enough.
	\end{enumerate}}
\end{lem}

 \begin{rem}
Observe that condition $(*_K)$ is very natural. Indeed, it is satisfied by any action having  no global fixed points and a free orbit. This is always the case for dynamical realizations of left orders on countable groups.
  
 \end{rem}
 
\noindent {\bf Dependency structure of results:} Theorem \ref{teo amalgama} is derived from Lemma \ref{lem amalgama}. We use Lemmas \ref{lem commutator semicont} and \ref{lem commutator} to deduce two auxiliary lemmas (\ref{lem taming} and \ref{lem removing}) that yield Theorems \ref{teo bonatti}, \ref{teo conjugacy} and \ref{teo dense orbit}. Theorem \ref{teo manija} is deduced from Lemma \ref{lem commutator}.


\subsection{Flexibility in amalgamated free products}
\label{sec proof teo amalgama}

In this section we show Theorem \ref{teo amalgama}.

Let $G=\F_n*_{w_1=w_2}\F_m=\langle x_1,...,x_{n+m}|w_1(x_1,...,x_n)=w_2(x_{n+1},...,y_{n+m})\rangle$, $n\geq 2$, $m\geq 1$. We can assume $w_1$ and $w_2$ are cyclically reduced, since conjugate words yield isomorphic amalgamated products. Take $\rho:G\to Homeo_{+}(\R)$ a representation without global fixed points and $p\in\R$. We will construct $\rho'$ another representation with no global fixed points so that $\rho'(x_i)$ and $\rho(x_i)$ coincide over $(-\infty,p]$ for each $i=1,...,n+m$, and such that $\rho'(w_1)$ and $\rho(w_1)$ are not semi-conjugated. Therefore we will get $\rho'$ a perturbation of $\rho$ not semi-conjugated to it.

Given $\rho\in Rep(G, Homeo_+(\R))$, define $\rho_1\in Rep({\F_n,Homeo_{+}(\R)})$ and $\rho_2\in Rep({\F_m,Homeo_{+}(\R)})$ as the restrictions of $\rho$ to the first and second factors of the amalgam decomposition of $G$. 

{\bf Case I:} $m\geq 2$.

\vs 

Take $q_0>q$ as the maximum of the $q_0$ given by Lemma \ref{lem amalgama} for the representations $\rho_1$ and $\rho_2$, and the point $p$. We will distinguish two cases. 

{\bf Subcase Ia:} $Fix(\rho(w_1))$ does not contain a neighbourhood of $+\infty$. 

\vs
Consider $h\in Homeo_{+}(\R)$ and $q>q_0$ such that:
\begin{itemize}

\item $\rho(w_1)(q)\neq q$.
\item $h$ coincides with $\rho_1(w_1)$ on $(-\infty,q]$.

\item $Fix(h)\cap(q,+\infty)$ consists of at most one point, and in case  $h$ fixes some point in $(q,\infty)$, then  $h$ transverses the diagonal at that point. We also impose that
$$ \left\{ \begin{array}{ll} \mbox{If $\rho(w_1)(x)>x$ for $x$ big enough} &  \mbox{$h(x)<x$ for $x$ big enough.} \\  \mbox{If $\rho(w_1)(x)<x$ for $x$ big enough} & \mbox{$h(x)>x$ for $x$ big enough.}  \\ \mbox{If $Fix(\rho(w_1))$ accumulates at $+\infty$} & \mbox{no further condition on $h$.} \end{array} \right.$$
\end{itemize}

Observe that it is possible to choose such $q$ since $\rho(w_1)$ has arbitrarily big points that are not fixed.

Since $n,m\geq 2$ we can apply Lemma \ref{lem amalgama} and obtain $\rho_1^*$ and $\rho_2^*$ such that  
\begin{itemize}
\item $\rho_1^{*}(x_i)$ coincide with $\rho_1(x_i)$ over $(-\infty,p]$ for $i=1,...,m$
\item $\rho_2^{*}(x_i)$ coincide with $\rho_2(x_i)$ over $(-\infty,p]$ for $i=m+1,...,m+n$
\item $\rho_1^*(w_1(x_1,...,x_n))=\rho_2^*(w_2(x_{n+1},...,x_{n+m}))=h$
\item $Fix(\rho_1^*)\subseteq Fix(\rho_1)$ and $Fix(\rho_2^*)\subseteq Fix(\rho_2)$
\end{itemize}

Define $\rho'\in Rep(G,Homeo_{+}(\R))$ as $\rho'(x_i)=\rho_1^*(x_i)$ for $i=1,...,m$ and $\rho'(x_i)=\rho_2^*(x_i)$ for $i=n+1,...,n+m$. Then $\rho'$ satisfies the thesis of Theorem \ref{teo amalgama}.  $\hfill\Box$

\vs

{\bf Subcase Ib:} $Fix(\rho(w_1))$ contains a neighbourhood of $+\infty$.

\vs

Assume first that either $\rho_1$ or $\rho_2$ has global fixed points accumulating on $+\infty$. Let's say $\rho_1$ does. Take $q\in Fix(\rho_1)$ with $q>p$ and $\rho_1(w_1)(x)=x$ for $x>q$. We can define $\rho'_1$ that agrees with $\rho_1$ over $(-\infty,q]$, and on $(q,\infty)$ we put an action without global fixed points but such that $w_1$ acts trivially. This can be done, for instance, by first sending  $\F_n$ to  an infinite cyclic homomorphic image where $w_1$ becomes trivial. In this way, the resulting action $\rho_1'$ is certainly not semi-conjugate to the initial $\rho_1$.

If on the other hand  both $\rho_1$ and $\rho_2$ have no global fixed points on a neighbourhood of $+\infty$, we use the Alexander trick (see Remark \ref{rem alexander}) on one factor. Concretely,
Take $q>p+1$ such that $\rho_1(w_1)(x)=x$ for $x>q-1$ and $Fix(\rho_2)\cap(q-1,+\infty)=\emptyset$. Consider $\phi:(-\infty,q)\to\R$ and orientation preserving homeomorphism that restricts to the identity over $(-\infty,q-1]$. Define $\rho'_1$ as $\phi\circ\rho_1\circ\phi^{-1}$ on $(-\infty,q]$ and as the trivial action on $[q,+\infty)$.

In both cases, we have $\rho'_1(w_1)=\rho_1(w_1)$, and thus can define $\rho'\in Rep(G,Homeo_{+}(\R))$ by exchanging $\rho_1$ for $\rho'_1$ (leaving $\rho_2$ as it is). This gives a representation that is not semi-conjugated to $\rho$, since we changed the behaviour near $+\infty$ of the global fixed points of the first factor. By construction we have $Fix(\rho')=\emptyset$ in both cases. $\hfill\Box$
\vs

{\bf Case II:} $m=1$. 

\vs
Now we have $w_2 = x_{n+1}^k$ for some $k\geq 1$. Notice that $Fix(\rho)=Fix(\rho_1)$ in this case. We take $\rho_1^*\in Rep({\F_n,Homeo_{+}(\R)})$ a representation that is not semi-conjugate to $\rho_1$, but satisfies that $\rho_1^*(x_i)$ and $\rho_1^*(w_1)$ agrees with $\rho(x_i)$ and $\rho(w_1)$ over $(-\infty,p]$ respectively. This can certainly be constructed by perturbing $\rho$ very  close to infinity. We can also demand that $\rho_1^*$ has no global fixed points on $(p,+\infty)$ (see \cite{ghys}). Now let $f$ be a $k$-th root of $\rho_1^*(w_1)$ that agrees with $\rho(x_{n+1})$ on $(-\infty,p]$. We let $\rho' \in Rep(G,Homeo_{+}(\R)) $ be defined as $\rho'(x_i)=\rho_1^*(x_i)$ for $i=1,...,n$ and $\rho'(x_{n+1})=f$. The representation $\rho'$ satisfies the thesis of Theorem \ref{teo amalgama}. $\hfill\Box$

\subsection{Orderings on groups with `handle' decomposition}
\label{sec  teo manija}

As announced in the Introduction, in this section we show

\begin{thm}\label{teo manija} {\em  Suppose $G$ is a countable left orderable group admitting a decomposition of the form $H*_{h=[a,b]}\F_2(a,b)$, with $h\neq id$. Then the space of left order of $G$ has no isolated points. }
\end{thm}

\noindent{\bf Proof:} Suppose $\preceq$ is a left order on a countable group $G$ admitting a decomposition of the form $H *_{h=[a,b]}\F_2(a,b)$, with $h\neq id$. To show that $\preceq$ is non-isolated, we need to show that given any finite set $F\subset G$, there is a left order $\preceq'$ which is different to $\preceq$ but such that,  for every $s\in F$, we have that $id \preceq' s$ iff $id\preceq s$.
\vsp

Fix $S$ a (perhaps infinite) generating set of $G$. We say that a finite set $F\subseteq G$ is {\em closed under prefix} if for every $g\in F$, there is a way to write $g$ as a product $s_1 \ldots s_n$, where each $s_j\in S$, and such that for every $i\leq n$, the element $s_i\ldots s_i$ also belongs to $F$. Clearly every finite  set admits an over set which is finite and closed under prefix (for instance, when $S$ is finite it suffices to take a ball containing $F$).

\vsp

Let $\rho$ be the dynamical realization of $\preceq$, and $F$ be a finite subset of $G$ which is closed under prefix. Let $K$ be a compact interval containing $s(0)$ for every $s\in F$.  Since $\rho$ is the dynamical realization of $\preceq$ we have that $id\prec g$ if and only if $0< \rho(g)(0)$. Since by construction $\rho$ has no global fixed points and the orbit of $0$ is free (that is, $0$ has trivial stabilizer),  we have that $\{x\in\R:\rho(h)(x)\neq x\}$ accumulates on $+\infty$ and on $-\infty$. In particular, condition $(*_K)$ holds for the pair $(\rho(a),\rho(b)$). We can then let $\tilde{f}$ and $\tilde{g}$ be the homeomorphisms provided by the conclusion of Lemma \ref{lem commutator} when applied to $(\rho(a),\rho(b))$ on the compact set $K$. With it we define $\rho':G\to Homeo_+(\R)$ by $\rho'(a)=\tilde{f}$, $\rho'(b)=\tilde{g}$ and $\rho'(g)=\rho(g)$ for every $g\in H$. Since $(\tilde{f},\tilde{g})\in\mathcal V_{\rho(h)}$, $\rho':G\to Homeo_+(\R)$ is a well defined homomorphism.

\vsp

Denote by $\preceq'$ the left order on $G$ defined by $$g_1\preceq' g_2 \text{ if and only if } \rho'(g_1)(0)\leq \rho'(g_2)(0). \footnote{As defined, the order $\preceq'$ is really a partial left invariant order, but since $Stab_{\rho'(G)}(0)$ is left orderable, we can extend the partial order $\preceq'$ to a total left order. See for instance \cite[\S2.1]{GOD}.} $$
Since $F$ is prefix closed, it follows from the definition of $\rho'$ that $\rho'(s)(0)=\rho(s)(0)$ for every $s\in F$. In particular, for every $s\in F$ we have that $id \prec' s$ if and only if $id\prec s$. So $\preceq'$ is close to $\preceq$.

Further, since Lemma \ref{lem commutator} ensures that the set of fixed points of $\langle \tilde{f},\tilde{g}\rangle$ is contained in the set of fixed points of $\langle \rho(a),\rho(b)\rangle$, we have that $\rho'(G)$ has no global fixed points. Now we check that $\preceq'$ and $\preceq$ are different. Indeed, if we let $w_n\in G$ be a sequence of elements such that $\rho'(w_n)(0)\geq n$, then for every $n$ large enough we have that $\rho'(w^{-1}_n b w_n)(0)=\rho'(w_n^{-1})\circ \tilde{g} \circ\rho'(w_n)(0)$ is either positive or negative depending on whether $\tilde{g}(x)-x$ is positive or negative for every $x$ large enough (see the conclusion of Lemma \ref{lem commutator}). This shows that, if we choose the appropriate behaviour of $\tilde{g}$ at infinity, then $\preceq'$ is different from $\preceq$, and therefore that $\preceq$ is non-isolated. $\hfill\square$

\subsection{Removing fixed points}
\label{sec proof teo bonatti}

In this section we give the proof of Theorem \ref{teo bonatti}. Throughout this section, $M$ will denote either the real line or the circle. As customary, for $a,b\in S^1$, the {\em interval} $(a,b)$ is the set of points $p\in S^1$ such that $(a, p , b)$ is clockwise oriented.

Let $q\geq 2$, and $\Sigma_q$ be a genus $q$ closed (orientable) surface. Our prefered presentation for $\pi_1(\Sigma_q)$ we will be $\langle a_1,b_1,\ldots ,a_q,b_q\mid [a_1,b_1]=w_1\rangle$, where $w_1=\Pi_{j\neq 1}[a_j,b_j]$. This is  the presentation induced from the amalgam decomposition
$$\pi_1(\Sigma_q)\simeq\F_2*_{[a_1,b_1]=w_1}\F_{2(q-1)},$$ where $\F_n$ is the free group of rank $n$. We will need the following definitions:

\begin{defi} Let $p\in M$ and $\phi\in Homeo_+(M)$. We will say that fixed point $p$ of $\phi$ is  of hyperbolic type if there exists a neighbourhood $V$ of $p$ such that $V- \{p\}$ has two connected components, and such that either $\phi^n$ or $\phi^{-n}$  shrinks $V$ to $\{p\}$ as $n\to \infty$. 
\end{defi}

In dynamical terms, that is to say that $p$ is either an {\em attracting} or {\em repelling} fixed point of $\phi$.

\begin{defi} Let $\rho\in Rep(\pi_1(\Sigma_q),Homeo_+(M))$ and $p\in Fix(\rho)$. We say that $\rho$ is {\em tame} on $p$ if $p$ is both an isolated fixed point of $\rho([a_1,b_1])$ and a fixed point of hyperbolic type for $\rho(b_1)$. In this case, if $V$ is a convex neighbourhood of $p$ with $Fix(\rho([a_1,b_1]))\cap V= Fix(\rho(b_1))\cap V=\{p\}$, we say that $\rho$ is tame on $p$ over $V$.

\end{defi}

The skeleton of the proof is the following: First we will show that any representation of $\pi_1(\Sigma_q)$ on $Homeo_{+}(M)$ can be approximated by representations whose global fixed points are isolated. Next we will approximate a representation with isolated global fixed points by one that is tame on each of them. Finally we show how to remove tame global fixed points by small perturbations.

\vs

\begin{defi} Let $\mathcal{F}$ be a family of closed intervals. We will say that $\mathcal{F}$ is locally finite if given a compact set $K$, only finitely many intervals in $\mathcal{F}$ intersect $K$.
\end{defi}

\begin{lem} \label{lem isolating} {\bf (Isolating):}  {\em Given $\rho\in Rep(\pi_1(\Sigma_q),Heomeo_+(M))$, and $U$ a neighbourhood of $\rho$ in the compact open topology, there exists $\rho'\in U$ such that $Fix(\rho')$ consists of isolated points.}
\end{lem}

\vsp

\noindent{\bf Proof:} Let $\rho\in Rep(\pi_1(\Sigma_q),Homeo_+(M)$. Given $\epsilon>0$, consider a locally finite family $\mathcal{F}$ of closed intervals of diameter less than $\epsilon$ that satisfy
\begin{itemize}
\item If $I,J\in \mathcal{F}$ are different then its interiors are disjoint.
\item If $x\in Fix(\rho)$ then there exists $I\in \mathcal{F}$ that contains $x$.
\item The endpoints of every $I\in\mathcal{F}$ are global fixed points of $\rho$.
\end{itemize}

Now, for each $I \in\mathcal{F}$ consider $\rho_{I}\in Rep(\pi_1(\Sigma_q),Homeo_+(I))$ a representation without global fixed points (except for the endpoints of $I$) and define $\rho'\in Rep(\pi_1(\Sigma_q),Homeo_+(M))$ as $\rho'(g)(x)=\rho(g)(x)$ if $x\notin \cup_{I\in \mathcal{F}}I$, and $\rho'(g)(x)=\rho_I(g)(x)$ if $x\in I$ for some $I\in \mathcal{F}$. Note that if $\epsilon$ is sufficiently small then $\rho'\in U$.

Finally, the local-finiteness of $\mathcal{F}$ implies that $Fix(\rho')$ is a  discrete set. $\hfill\square$

\vs

\begin{lem} \label{lem taming} {\bf (Taming):} {\em Let $\rho\in Rep(\pi_1(\Sigma_q),Homeo_+(M))$, $p\in Fix(\rho)$ and $V$ a neighbourhood of $p$ in $M$ such that $Fix(\rho)\cap V=\{p\}$. Then there exists $\rho'\in Rep(\pi_1(\Sigma_q),Homeo_+(M))$ such that
\begin{itemize}
\item $\rho'(a_i)_{|V^c}=\rho(a_i)_{|V^c}$ and $\rho'(b_i)_{|V^c}=\rho(b_i)_{|V^c}$ for $i=1,...,q$, where $V^c$ denotes the complement of $V$

\item $\rho'$ is tame on $p$
\item $Fix(\rho')$ is included in $Fix(\rho)$.
\end{itemize}
}
\end{lem}

\noindent{\bf Proof:} Consider $(\alpha,p)$ a subset of $M$ homeomorphic to $\R$ where $\alpha=-\infty$ or $\alpha\in Fix(\rho)$, and such that $Fix(\rho)\cap (\alpha,p)=\emptyset$. Take $\Phi\colon\R\to (\alpha,p)$ an orientation preserving homeomorphism. We will consider the representation $\theta\in Rep(\pi_1(\Sigma_q),Homeo_+(\R))$ defined as $\theta(g)(x)=\Phi^{-1}\rho(g)\Phi(x)$. This is equivalent to consider $\rho$ acting on the interval $(\alpha,p)$.

Take $k\in\R$ such that $\Phi([k,+\infty))\subseteq V$. Apply Lemma \ref{lem commutator semicont} to get a perturbation $\theta_1$ such that $\theta_1(a_i)_{|(-\infty,k)}=\theta(a_i)_{|(-\infty,k)}$ and  $\theta_1(b_i)_{|(-\infty,k)}=\theta(b_i)_{|(-\infty,k)}$ for $i=1,\ldots,q$, and also that $\theta_1([a_1,b_1])(x)=x+1$ for $x$ big enough and $Fix(\theta_1)=\emptyset$.

Since $\theta_1([a_1,b_1])(x)=x+1$ for $x$ large enough, the pair $(\theta_1(a_1),\theta_1(b_1))$ satisfies condition $(*_{(-\infty,k]})$. Therefore we can apply Lemma \ref{lem commutator} and find $f_1$ and $g_1$ perturbations of $\theta_1(a_1)$ and $\theta_1(b_1)$ supported on $(k,+\infty)$, such that $[f_1,g_1]=\theta_1([a_1,b_1])$, $Fix(\langle f_1,g_1\rangle)\subseteq Fix(\langle\theta_1(a_1),\theta_1(b_1)\rangle)$ and $g_1(x)>x$ for $x$ large enough.

Since $[f_1,g_1]=\theta_1([a_1,b_1])$ we can define a representation $\theta'\in Rep(\pi_1(\Sigma_q),Homeo_+(\R))$ by $\theta'(a_1)=f_1$, $\theta'(b_1)=g_1$, $\theta'(a_i)=\theta_1(a_i)$ and $\theta'(b_i)=\theta_1(b_i)$ for $i=2,\ldots,q$. By construction, we have that there is a neighbourhood of $+\infty$ where $\theta'(b_1)=g_1$ is increasing and $\theta'([a_1,b_1])=\theta_1([a_1,b_1])$ has no fixed points.

Finally, since $Fix(\langle f_1,g_1\rangle)\subseteq Fix(\langle\theta_1(a_1),\theta_1(b_1)\rangle)$ we have $$Fix(\theta')=Fix(\langle f_1,g_1\rangle)\cap Fix(\langle \theta_1(a_2),\theta_1(b_2),\ldots,\theta_1(a_q),\theta_1(b_q)\rangle) \subseteq Fix(\theta_1)=\emptyset.$$

We define $\bar{\rho}$ by $\bar{\rho}(g)(x)=\rho(g)(x)$ if $x\notin (\alpha,p)$ and $\bar{\rho}(g)(x)=\Phi\theta'(g)\Phi^{-1}(x)$ if $x\in (\alpha,p)$. Note that $\bar{\rho}(a_i)_{|V^c}=\rho(a_i)_{|V^c}$ and $\bar{\rho}(b_i)_{|V^c}=\rho(b_i)_{|V^c}$ for $i=1,...,q$. By this construction we have that $Fix(\bar{\rho})\subseteq Fix(\rho)$, and also that $p$ has a neighbourhood $V'$ so that in the left component of $V'-\{p\}$ there are no fixed points of $\bar{\rho}([a_1,b_1])$ and $\bar{\rho}(b_1)$ is increasing.

We repeat the same procedure on the other side of $p$ and get $\rho'$. We will have that $p$ is an isolated point of $Fix(\rho'([a_1,b_1]))$. Moreover, taking the right choice when applying Lemma \ref{lem commutator}, we will have that $p$ is of hyperbolic type for $\rho'(b_1)$ and so $\rho'$ is tame on $p$. Finally, note that $Fix(\rho')\subseteq Fix(\bar\rho)\subseteq Fix(\rho)$. $\hfill\square$

\vs

The following Lemma shows how to remove tame global fixed points

\vs
\begin{lem} \label{lem removing} {\bf (Removing):}  {\em Let $\rho\in Rep(\pi_1(\Sigma_q),Homeo_+(M))$ and $p\in Fix(\rho)$ that is tamed over an interval $V$. Then we can construct $\bar{\rho}\in Rep(\pi_1(\Sigma_q),Homeo_+(M))$ such that \begin{itemize}
\item $\bar{\rho}(a_i)_{|V^c}=\rho(a_i)_{|V^c}$ and $\bar{\rho}(b_i)_{|V^c}=\rho(b_i)_{|V^c}$ for $i=1,...,q$
\item $Fix(\bar{\rho})\cap V=\emptyset$
\end{itemize}}
\end{lem}
{\bf Proof:}
Since $\rho$ is tame over $V$ we can construct $g$ a perturbation of $\rho(b_1)$ supported on $V$ such that the graph of $g_{|V}$ transverses the graphs of the identity and of $\rho([a_1,b_1]^{-1})$ only once, at different points $p_1$ and $p_2$ respectively. These points are then the only fixed points of $g$ and $\rho([a_1,b_1])g$ on the interval $V$. In particular $Fix(g)\cap Fix(\rho([a_1,b_1]))\cap V=\emptyset$.

Notice that a homeomorphism $\psi$ with $\psi_{|V^c} = \rho(a_1)_{|V^c}$ and $\psi(p_1)=p_2$ is a weak conjugation from $g$ to $\rho([a_1,b_1])g$ that is strong on each component of $V^c$. Since these components are separated by a fixed point of $g$, the arguments for Lemma \ref{promotion} also work in this case. So we get $f\in Homeo_{+}(\R)$ such that ${f}_{|V^c}=\rho(a_1)_{|V^c}$ and that conjugates $g$ to $\rho([a_1,b_1])g$.

Finally, we define $\bar{\rho}$ as
\begin{itemize}
\item $\bar{\rho}(a_i)=\rho(a_i)$ and $\bar{\rho}(b_i)=\rho(b_i)$ for $i=2,...,q$
\item $\bar{\rho}(a_1)=f$ and $\bar{\rho}(b_1)=g$
\end{itemize}
$\hfill\square$

Now we are in position to finish the
\vs

{\bf Proof of Theorem \ref{teo bonatti}:} Let $\rho\in Rep(\pi_1(\Sigma_q),Homeo_+(M))$ and $U$ a neighbourhood of $\rho$ in the compact open topology. First we apply the Isolating Lemma to find $\rho_1\in U$ with isolated global fixed points.

For each $p\in Fix(\rho_1)$ take a neighbourhood $V_p$ so that $p$ is the only global fixed point of $\rho_1$ on it. We can also assume they are pairwise disjoint. On each $V_p$ we apply the perturbation of the Taming Lemma followed by that of the Removing Lemma. We can do this recursively (for some order of $Fix(\rho_1)$) and take the limit. This will be the representation $\rho'$ in the statement of Theorem \ref{teo bonatti}. It is clear that $Fix(\rho')=\emptyset$. Finally, notice that by taking the $V_p$ small enough we can guarantee that $\rho'\in U$.  $\hfill\square$

\subsection{Construction of a dense orbit  in $Rep_\#(\pi_1(\Sigma),Homeo_+(\R))$}

In this section we give the proof of Theorem \ref{teo conjugacy}.

Fix an orientation preserving homeomorphism $\Phi:\R\to (0,1)$, and for each $n\in\Z$ let $\Phi_n(x)=\Phi(x)+n$.

We will write $Rep_{\#}:=Rep_{\#}(\pi_1(\Sigma_q),Homeo_{+}(\R))$. Notice that $Rep_{\#}$ is 
separable (by Remark \ref{separability}), so we can consider $Q\subseteq Rep_{\#}$ a dense countable subset. Let $\{\rho_n:n\in\Z\}$ be a sequence in $Q$ that repeats every element infinitely often. We will define $\theta_0\in Rep(\pi_1(\Sigma_q),Homeo_{+}(\R))$ as follows:
\begin{itemize}
\item Each $n\in\Z$ is a global fixed point.
\item On the interval $(n,n+1)$, define $\theta_0(g) = \Phi_n \rho_n(g) \Phi_n^{-1}$ for all $g\in \pi_1(\Sigma_q)$.
\end{itemize}

Since $\rho_n\in Rep_{\#}$, we see that $Fix(\theta_0)=\Z$.

For each $n\in\Z$ we consider the interval $V_n = (n-2^{-|n|-1},n+2^{-|n|-1})$. Then $V_n$ is a convex neighbourhood of $n$ with $diam V_n < 2^{-|n|}$, disjoint with any other $V_m$. As in the proof of Theorem \ref{teo bonatti}, we apply the Taming Lemma followed by the Removing Lemma on each $V_n$, obtaining a representation $\theta$ without global fixed points.

We claim that the conjugacy class of $\theta$ is dense in $Rep_{\#}$. To show this, it is enough to prove that for every $\rho\in Q$ and every $m>0$ there is a conjugate $\bar \theta$ of $\theta$ so that $\bar \theta(a_i)_{|[-m,m]} = \rho(a_i)_{|[-m,m]}$ and $\bar \theta(b_i)_{|[-m,m]} = \rho(b_i)_{|[-m,m]}$ for $i=1,\ldots,q$.

Take $n\in\Z$ so that $\rho_n=\rho$ and $\Phi_n([-m,m])$ is disjoint with $V_n\cup V_{n+1}$. This is possible since the sequence $\{\rho_n:n\in\Z\}$ repeats $\rho$ infinitely many times, and $diam V_n$ goes to $0$ as $|n|\to+\infty$. Take $\psi\in Homeo_+(\R)$ that agrees with $\Phi_n$ on $[-m,m]$. Then $\psi^{-1}\theta(a_i)\psi_{|[-m,m]}=\rho(a_i)_{|[-m,m]}$ and $\psi^{-1}\theta(b_i)\psi_{|[-m,m]}=\rho(b_i)_{|[-m,m]}$ for $i=1,\ldots,q$.

Finally, applying Theorem \ref{teo bonatti} we get that $Rep_{\#}$ is dense in $Rep$ and therefore the conjugacy class of $\theta$ is dense in $Rep(\pi_1(\Sigma), Homeo_{+}(\R))$. $\mbox{ }$ $\hfill\square$

\subsection{Construction of a dense orbit in $\mathcal{LO}(\pi_1(\Sigma))$}

In this section we give the proof of Theorem \ref{teo dense orbit}. The construction follows closely the one for Theorem \ref{teo conjugacy}. 

We take $Q$ a countable dense subset of $\mathcal{LO}(\pi_1(\Sigma))$ (this certainly exists since the space of left orders of countable groups is compact and metrizable,therefore separable,  see \S \ref{sec orderings}), and $\{\rho_n:n\in\Z-\{0\}\}$ a sequence of dynamical realizations of the orders in $Q$, repeating each representation infinitely often. Let $\rho_0$ be a representation of $\pi_1(\Sigma)$ by translations with  dense orbits (e.g. translations by lengths that are linearly independent over $\mathbb{Q}$).

Again let $\Phi:\R\to (-1,1)$ be an orientation preserving homeomorphism with $\Phi(0)=0$, and for each $n\in\Z$ let $\Phi_n(x)=\Phi(x)+2n$. We define the representation $\theta_0$ as follows:
\begin{itemize}
\item Each odd integer is a global fixed point.
\item On the interval $(2n-1,2n+1)$, define $\theta_0(g) = \Phi_n \rho_n(g) \Phi_n^{-1}$ for all $g\in \pi_1(\Sigma_q)$.
\end{itemize}

For each odd integer $n$ we take a convex neighbourhood $V_n$ with $diam V_n < 2^{-|n|}$, and we use the Taming and Removing Lemmas as in the proof of Theorem \ref{teo bonatti} to remove the global fixed points with a perturbation supported on the $V_n$. Let $\theta$ be the representation thus obtained.

We will check that if $V_1$ and $V_{-1}$ are small enough, then the orbit of $0$ under $\theta$ is dense.

Let $S$ be a generating set of $\pi_1(\Sigma)$. For $M>0$ and a representation $\rho$ consider the {\em local orbit} of $0$ on $[-M,M]$, that is the set $L_M(\rho)$ of points of the form $\rho(g)(0)$ where $g=s_1\cdots s_k$ with $s_j\in S$ and so that $\rho(s_i\cdots s_k)(0)\in [-M,M]$ for all $i=1,\ldots,k$. Notice that a perturbation of $\rho$ outside $[-M,M]$ does not change the local orbit $L_M$.

By our choice of $\rho_0$, we have that for $M$ big enough the closure of $L_M(\rho_0)$ contains a neighbourhood of $0$. So by taking $V_1$ and $V_{-1}$ small enough (disjoint from $\Phi([-M,M])$) we get that there is a neighbourhood of $0$ contained in the closure of some local orbit of $0$ under $\theta$.  Therefore the closure of the orbit of $0$ under $\theta$ is both open and closed, so the orbit is dense.

Define $\prec$ by $$g_1\preceq g_2 \text{ if and only if } \theta(g_1)(0)\leq \theta(g_2)(0)$$
As in the proof of Theorem \ref{teo manija}, $\prec$ is really a partial left order, if $0$ has no trivial stabilizer, but it can be extended to a total left order. Any such extension will have a dense orbit under conjugation.

Let $\tilde{\prec}$ be any element of $\mathcal{LO}(\pi_1(\Sigma))$ and a finite subset $F$ of $\pi_1(\Sigma)$. We can assume $F$ is {\em closed under prefix} (see the beginin of  \S \ref{sec teo manija} for definition). Take  $\prec'\in Q$ an element that agrees with $\tilde{\prec}$ on $F$, and let $\rho'$ be its dynamical realization. Consider $K$ the convex closure of the $\rho'(f)(0)$ for $f\in F$. Since there are infinitely many repetitions of $\rho'$ in the sequence $\{\rho_n:n\in\Z-\{0\}\}$, we can take one with $n$ big enough so that $\Phi_n(K)$ is disjoint with $V_{2n-1}\cup V_{2n+1}$. So we get that $$ f_1\prec' f_2 \text{ if and only if } \theta(f_1)(2n)<\theta(f_2)(2n) \text{ for } f_1,f_2\in F$$ noticing that $2n = \Phi_n(0)$.

But now, since $F$ is finite and $\theta(f)(2n)\neq 2n$ for every $f\in F$, we have that there is a neighbourhood $U$ of $2n$ so that, for any $p\in U$ and $f_1,f_2\in F$, we have that $f_1\prec' f_2$ if and only if $\theta(f_1)(p)<\theta(f_2)(p)$. Since by construction the orbit of $0$ under $\theta$ is dense, we can take $g\in\pi_1(\Sigma)$ with $\theta(g)^{-1}(0)\in U$. It follows that $\prec_g$ agrees with $\prec'$, and therefore with $\tilde{\prec}$, on $F$, as desired. $\hfill\square$

\section{Proof of Lemma \ref{lem amalgama}}
\label{sec proof lem amalgama}

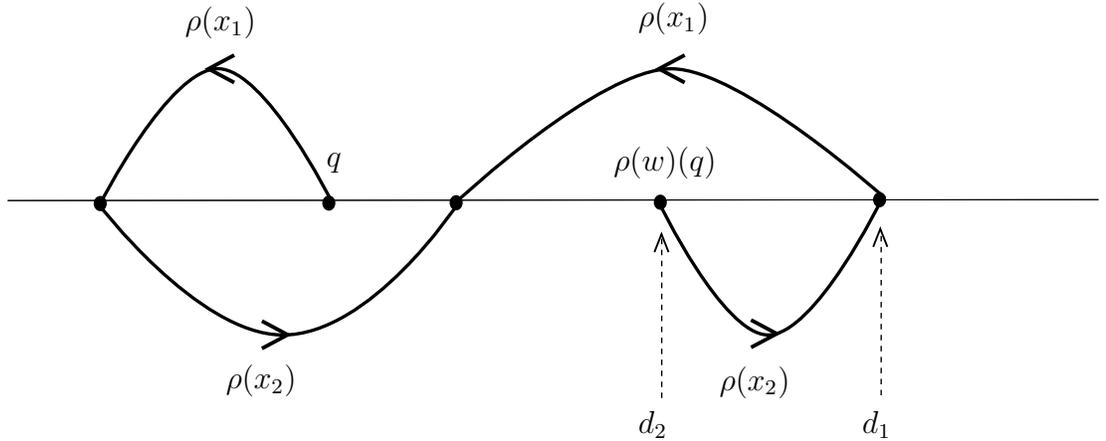
\begin{figure}[h] 

\begin{tikzpicture}[y=0.80pt, x=0.75pt, yscale=-1.000000, xscale=1.000000, inner sep=0pt, outer sep=0pt]
\draw (-200,0) node  {$ $};

  \path[draw=black,fill=black,line join=miter,line cap=butt,even odd rule,line
    width=0.382pt] (-171.2731,75.1647) -- (378.8367,74.7268);
  \path[draw=black,line join=miter,line cap=butt,miter limit=4.00,line
    width=1.255pt] (-8.0076,75.4398) .. controls (-8.0076,75.4398) and
    (-41.2484,13.0985) .. (-65.0397,12.7615) .. controls (-89.8107,12.4106) and
    (-124.8951,77.0509) .. (-124.8951,77.0509);
  \path[draw=black,line join=miter,line cap=butt,miter limit=4.00,line
    width=1.255pt] (271.1982,73.9901) .. controls (271.1982,73.9901) and
    (198.8783,13.1439) .. (163.7557,12.8069) .. controls (127.1868,12.4560) and
    (55.2448,75.6012) .. (55.2448,75.6012);
  \path[draw=black,line join=miter,line cap=butt,miter limit=4.00,line
    width=1.255pt] (-126.1844,76.0828) .. controls (-126.1844,76.0828) and
    (-75.6945,138.4241) .. (-32.8991,138.7611) .. controls (11.6586,139.1120) and
    (55.1281,77.4619) .. (55.1281,77.4619);
  \path[draw=black,line join=miter,line cap=butt,miter limit=4.00,line
    width=1.255pt] (157.1762,76.0828) .. controls (157.1762,76.0828) and
    (189.1463,138.4241) .. (212.0280,138.7611) .. controls (235.8521,139.1120) and
    (269.5952,74.4717) .. (269.5952,74.4717);
  \path[cm={{2.3,0.0,0.0,2.3,(-134.91631,-126.13677)}},draw=black,fill=black,miter
    limit=4.00,line width=0.546pt] (4.4793,88.1312) ellipse (0.0336cm and
    0.0398cm);
  \path[cm={{2.3,0.0,0.0,2.3,(-134.91631,-126.13677)}},draw=black,fill=black,miter
    limit=4.00,line width=0.546pt] (54.6094,88.0378) ellipse (0.0336cm and
    0.0398cm);
  \path[cm={{2.3,0.0,0.0,2.3,(-134.91631,-126.13677)}},draw=black,fill=black,miter
    limit=4.00,line width=0.546pt] (82.5265,87.9670) ellipse (0.0336cm and
    0.0398cm);
  \path[cm={{2.3,0.0,0.0,2.3,(-134.91631,-126.13677)}},draw=black,fill=black,miter
    limit=4.00,line width=0.546pt] (127.2976,87.8641) ellipse (0.0336cm and
    0.0398cm);
  \path[cm={{2.3,0.0,0.0,2.3,(-134.91631,-126.13677)}},draw=black,fill=black,miter
    limit=4.00,line width=0.546pt] (175.3797,87.3492) ellipse (0.0336cm and
    0.0398cm);
  \path[xscale=0.919,yscale=1.088,fill=black] (-88.8867,-1.6066) node[above right]
    (text10353) {$\rho(x_1)$};
  \path[xscale=0.919,yscale=1.088,fill=black] (-66.6052,154.1979) node[above
    right] (text10363) {$\rho(x_2)$};
  \path[xscale=0.919,yscale=1.088,fill=black] (282.6144,172.1892) node[above
    right] (text10377) {$d_1$};
  \path[xscale=0.919,yscale=1.088,fill=black] (159.7385,172.1892) node[above
    right] (text10381) {$d_2$};
  \path[fill=black] (79.2904,194.5855) node[above right] (text10385) {};
  \path[draw=black,dash pattern=on 2.08pt off 2.08pt,line join=miter,line
    cap=butt,miter limit=4.00,line width=0.519pt] (158.5771,168.6210) --
    (158.5771,90.7277);
  \path[draw=black,dash pattern=on 2.08pt off 2.08pt,line join=miter,line
    cap=butt,miter limit=4.00,line width=0.519pt] (269.2472,167.3109) --
    (269.2472,89.4176);
  \path[xscale=0.919,yscale=1.088,fill=black] (-11.3299,56.8678) node[above right]
    (text7980) {$q$};
  \path[xscale=0.886,yscale=1.128,fill=black] (151.7320,57.4052) node[above right]
    (text7988) {$\rho(w)(q)$};
  \path[color=black,draw=black,line join=miter,line cap=butt,miter limit=4.00,line
    width=1.439pt] (-43.0527,132.1788) -- (-30.0705,138.6700) --
    (-43.0527,145.1610);
  \path[color=black,draw=black,line join=miter,line cap=butt,miter limit=4.00,line
    width=1.439pt] (-57.0230,6.3432) -- (-70.0052,12.8343) -- (-57.0230,19.3254);
  \path[color=black,draw=black,line join=miter,line cap=butt,miter limit=4.00,line
    width=1.439pt] (170.1659,6.3432) -- (157.1837,12.8343) -- (170.1659,19.3254);
  \path[color=black,draw=black,line join=miter,line cap=butt,miter limit=4.00,line
    width=1.439pt] (203.6214,131.6795) -- (216.6036,138.1706) --
    (203.6214,144.6617);
  \path[color=black,draw=black,line join=miter,line cap=butt,miter limit=4.00,line
    width=0.716pt] (272.3840,97.0758) -- (268.8086,88.5423) -- (265.2333,97.0758);
  \path[color=black,draw=black,line join=miter,line cap=butt,miter limit=4.00,line
    width=0.716pt] (162.0298,99.5725) -- (158.4544,91.0390) -- (154.8791,99.5725);
  \path[xscale=0.919,yscale=1.088,fill=black] (204.2922,154.8304) node[above
    right] (text3089) {$\rho(x_2)$};
  \path[xscale=0.919,yscale=1.088,fill=black] (159.7320,-3.3401) node[above right]
    (text3095) {$\rho(x_1)$};

\end{tikzpicture}

\caption{The picture shows a possible example of a sequence $S(\rho,w,q)$ for $w=x_2^{-1}x_1^{-1}x_2x_1$. We denote $d_i = d_w(S(\rho,w,q),i)$. } \label{orbita1}
\end{figure}

Let $w\in\F_n=\langle x_1,...,x_n\rangle$ cyclically reduced such that $w=a_m ... a_1$ with $a_i\in\{x_1^{\pm 1},...,x_n^{\pm n}\}$. We define $w_0=e$ and $w_j=a_{j}...a_{1}$ for $0<j\leq m$.
If $\rho\in Rep(\F_n,Homeo_{+}(\R))$ and $x\in\R$ we will be interested in the sequence $S(\rho,w,x)=(\rho(w_0)(x),...,\rho(w_m)(x))$. 
For each generator $x_i$ we will look at the minimum point from which we can perturb $\rho(x_i)$ without changing the sequence $S(\rho,w,x)$. With this in mind, for a general sequence $S=(s_0,...,s_m)$ and $k\in\{1,...,n\}$ we define $D_w(S,k)=\{s_j:a_{j+1}=k \ or \ a_{j}=x_k^{-1}\}$ and $d_w(S,k) = max D_w(S,k)$. (Figure \ref{orbita1} provides an example).

 Recall $p$ and $\rho$ from the statement of the Lemma \ref{lem amalgama}, and let $f_i=\rho(x_i)$ for $i=1,...,n$. We take $q_0$ such that $max\{\rho(u)(p): u\in B_1(\langle x_1,...,x_n \rangle) \}$ is less than every point in $S(\rho,w,q_0)$ and $S(\rho,w^{-1},q_0)$.

Take $q>q_0$ and $h$ an homeomorphism as in the statement of Lemma \ref{lem amalgama}. Let $d_i = d_w(S(\rho,w,q),i)$ (See figure \ref{orbita1}). We will first define $\overline{\rho}\in Rep(\F_n,Homeo_{+}(\R))$ such that $\overline{\rho}(w)$ is conjugated to $h$ and $\overline{\rho}(x_i) = g_i$ agrees with $f_i$ over $(-\infty,d_i]$ for $i=1,...,n$. We will do this by defining each $g_i$ on a discrete subset of $(d_i,+\infty)$ and then extend by interpolation.

\begin{lem} \label{lem interno} There exists $\overline{\rho}\in Rep(\F_n,Homeo_{+}(\R))$ such that:
\begin{itemize}
\item $\overline{\rho}(x_i)$ agrees with $\rho(x_i)$ over $(-\infty,d_i]$ for $i=1,...,n$.
\item $\overline{\rho}(w)$ is weakly conjugated to $h$, by a map that coincides with the identity on $(-\infty,q]$. Moreover, this weak conjugation is strong on $(-\infty,q]$.
\item $Fix(\overline{\rho})\subseteq Fix(\rho)$.
\end{itemize}
\end{lem}


{\bf Proof:}

{\bf Case Ia:} $h(q)>q$ and $Fix(h)\cap(q,+\infty)=\emptyset$.
\vs


 Take $r_1$ so that $q<r_1<h(q)$.

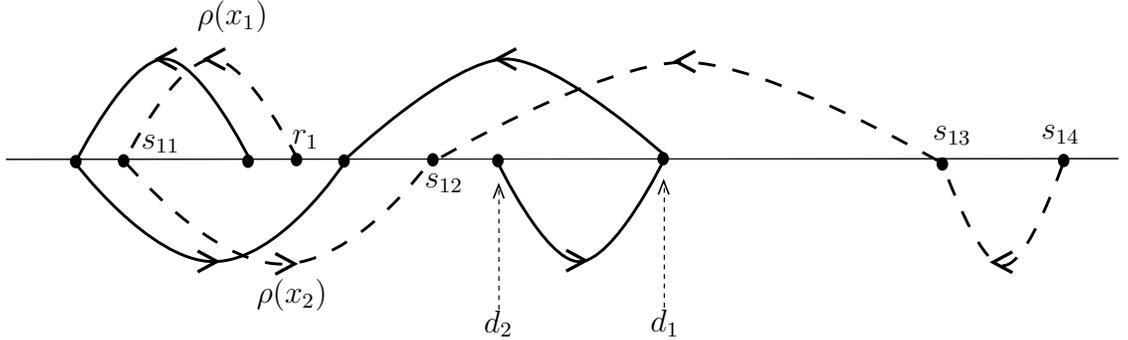
\begin{figure}[h]

\begin{tikzpicture}[y=0.70pt, x=0.650pt, yscale=-1.000000, xscale=1.000000, inner sep=0pt, outer sep=0pt]
\draw (-280,0) node  {$ $};

  \path[draw=black,fill=black,line join=miter,line cap=butt,even odd rule,line
    width=0.387pt] (-236.2435,80.2725) -- (411.0474,79.8917);
  \path[draw=black,line join=miter,line cap=butt,miter limit=4.00,line
    width=1.091pt] (-94.2735,80.5117) .. controls (-94.2735,80.5117) and
    (-123.1786,26.3019) .. (-143.8666,26.0088) .. controls (-165.4066,25.7037) and
    (-195.9148,81.9126) .. (-195.9148,81.9126);
  \path[draw=black,line join=miter,line cap=butt,miter limit=4.00,line
    width=1.091pt] (148.5142,79.2511) .. controls (148.5142,79.2511) and
    (85.6273,26.3413) .. (55.0859,26.0483) .. controls (23.2869,25.7431) and
    (-39.2714,80.6520) .. (-39.2714,80.6520);
  \path[draw=black,line join=miter,line cap=butt,miter limit=4.00,line
    width=1.091pt] (-197.0360,81.0708) .. controls (-197.0360,81.0708) and
    (-153.1317,135.2806) .. (-115.9183,135.5736) .. controls (-77.1724,135.8788)
    and (-39.3729,82.2700) .. (-39.3729,82.2700);
  \path[draw=black,line join=miter,line cap=butt,miter limit=4.00,line
    width=1.091pt] (49.3646,81.0708) .. controls (49.3646,81.0708) and
    (77.1647,135.2806) .. (97.0618,135.5737) .. controls (117.7784,135.8788) and
    (147.1203,79.6699) .. (147.1203,79.6699);
  \path[cm={{2.0,0.0,0.0,2.0,(-204.62889,-94.77233)}},draw=black,fill=black,miter
    limit=4.00,line width=0.546pt] (4.4793,88.1312) ellipse (0.0336cm and
    0.0398cm);
  \path[cm={{2.0,0.0,0.0,2.0,(-204.62889,-94.77233)}},draw=black,fill=black,miter
    limit=4.00,line width=0.546pt] (54.6094,88.0378) ellipse (0.0336cm and
    0.0398cm);
  \path[cm={{2.0,0.0,0.0,2.0,(-204.62889,-94.77233)}},draw=black,fill=black,miter
    limit=4.00,line width=0.546pt] (82.5265,87.9670) ellipse (0.0336cm and
    0.0398cm);
  \path[cm={{2.0,0.0,0.0,2.0,(-204.62889,-94.77233)}},draw=black,fill=black,miter
    limit=4.00,line width=0.546pt] (127.2976,87.8641) ellipse (0.0336cm and
    0.0398cm);
  \path[cm={{2.0,0.0,0.0,2.0,(-204.62889,-94.77233)}},draw=black,fill=black,miter
    limit=4.00,line width=0.546pt] (175.3797,87.3492) ellipse (0.0336cm and
    0.0398cm);
  \path[xscale=0.919,yscale=1.088,fill=black] (-135.8893,10.7175) node[above
    right] (text10353) {$\rho(x_1)$};
  \path[xscale=0.919,yscale=1.088,fill=black] (-98.2438,149.1704) node[above
    right] (text10363) {$\rho(x_2)$};
  \path[xscale=0.919,yscale=1.088,fill=black] (151.0381,161.6279) node[above
    right] (text10377) {$d_1$};
  \path[xscale=0.919,yscale=1.088,fill=black] (45.6915,161.9745) node[above right]
    (text10381) {$d_2$};
  \path[fill=black] (-18.3622,184.1166) node[above right] (text10385) {};
  \path[draw=black,dash pattern=on 1.81pt off 1.81pt,line join=miter,line
    cap=butt,miter limit=4.00,line width=0.452pt] (50.5828,161.5388) --
    (50.5828,93.8055);
  \path[draw=black,dash pattern=on 1.81pt off 1.81pt,line join=miter,line
    cap=butt,miter limit=4.00,line width=0.452pt] (146.8177,160.3996) --
    (146.8177,92.6662);
  \path[draw=black,dash pattern=on 6.55pt off 6.55pt,line join=miter,line
    cap=butt,miter limit=4.00,line width=1.091pt] (-66.3107,80.6969) .. controls
    (-66.3107,80.6969) and (-95.2157,26.4871) .. (-115.9038,26.1940) .. controls
    (-137.4438,25.8889) and (-167.9520,82.0978) .. (-167.9520,82.0978);
  \path[draw=black,dash pattern=on 6.55pt off 6.55pt,line join=miter,line
    cap=butt,miter limit=4.00,line width=1.091pt] (-166.2330,82.3915) .. controls
    (-166.2330,82.3915) and (-117.5322,136.6013) .. (-76.2534,136.8943) ..
    controls (-33.2746,137.1995) and (8.6544,83.5907) .. (8.6544,83.5907);
  \path[draw=black,dash pattern=on 6.55pt off 6.55pt,line join=miter,line
    cap=butt,miter limit=4.00,line width=1.091pt] (305.8542,80.9799) .. controls
    (305.8542,80.9799) and (207.1108,28.0702) .. (159.1557,27.7771) .. controls
    (109.2256,27.4720) and (10.9984,82.3808) .. (10.9984,82.3808);
  \path[draw=black,dash pattern=on 6.55pt off 6.55pt,line join=miter,line
    cap=butt,miter limit=4.00,line width=1.091pt] (378.5798,83.2010) .. controls
    (378.5798,83.2010) and (358.6166,137.4108) .. (344.3285,137.7038) .. controls
    (329.4520,138.0090) and (308.3816,81.8001) .. (308.3816,81.8001);
  \path[cm={{2.0,0.0,0.0,2.0,(-204.62889,-94.77233)}},draw=black,fill=black,miter
    limit=4.00,line width=0.546pt] (18.3892,87.9264) ellipse (0.0336cm and
    0.0398cm);
  \path[cm={{2.0,0.0,0.0,2.0,(-204.62889,-94.77233)}},draw=black,fill=black,miter
    limit=4.00,line width=0.546pt] (68.7008,87.5954) ellipse (0.0336cm and
    0.0398cm);
  \path[cm={{2.0,0.0,0.0,2.0,(-204.62889,-94.77233)}},draw=black,fill=black,miter
    limit=4.00,line width=0.546pt] (108.3823,87.6503) ellipse (0.0336cm and
    0.0398cm);
  \path[cm={{2.0,0.0,0.0,2.0,(-204.62889,-94.77233)}},draw=black,fill=black,miter
    limit=4.00,line width=0.546pt] (256.6302,88.7908) ellipse (0.0336cm and
    0.0398cm);
  \path[cm={{2.0,0.0,0.0,2.0,(-204.62889,-94.77233)}},draw=black,fill=black,miter
    limit=4.00,line width=0.546pt] (291.9195,87.8355) ellipse (0.0336cm and
    0.0398cm);
  \path[xscale=0.919,yscale=1.088,fill=black] (-75.8664,68.4258) node[above right]
    (text7980) {$r_1$};
  \path[xscale=0.919,yscale=1.088,fill=black] (-171.3392,70.6074) node[above
    right] (text7984) {$s_{11}$};
  \path[xscale=0.919,yscale=1.088,fill=black] (8.0606,90.3040) node[above right]
    (text7988) {$s_{12}$};
  \path[xscale=0.919,yscale=1.088,fill=black] (329.8645,66.2221) node[above right]
    (text7988-0) {$s_{13}$};
  \path[xscale=0.919,yscale=1.088,fill=black] (398.4135,66.0459) node[above right]
    (text7988-1) {$s_{14}$};
  \path[color=black,draw=black,line join=miter,line cap=butt,miter limit=4.00,line
    width=1.251pt] (-124.7475,129.8499) -- (-113.4586,135.4944) --
    (-124.7475,141.1388);
  \path[color=black,draw=black,line join=miter,line cap=butt,miter limit=4.00,line
    width=1.251pt] (-136.8956,20.4277) -- (-148.1844,26.0721) --
    (-136.8956,31.7166);
  \path[color=black,draw=black,line join=miter,line cap=butt,miter limit=4.00,line
    width=1.251pt] (-108.6733,20.4277) -- (-119.9622,26.0721) --
    (-108.6733,31.7166);
  \path[color=black,draw=black,line join=miter,line cap=butt,miter limit=4.00,line
    width=1.251pt] (60.6600,20.4277) -- (49.3711,26.0721) -- (60.6600,31.7166);
  \path[color=black,draw=black,line join=miter,line cap=butt,miter limit=4.00,line
    width=1.251pt] (164.8653,21.7303) -- (153.5764,27.3747) -- (164.8653,33.0192);
  \path[color=black,draw=black,line join=miter,line cap=butt,miter limit=4.00,line
    width=1.251pt] (349.4081,130.2842) -- (338.1192,135.9286) --
    (349.4081,141.5730);
  \path[color=black,draw=black,line join=miter,line cap=butt,miter limit=4.00,line
    width=1.251pt] (-79.5897,131.1526) -- (-68.3009,136.7970) --
    (-79.5897,142.4415);
  \path[color=black,draw=black,line join=miter,line cap=butt,miter limit=4.00,line
    width=1.251pt] (89.7517,129.4157) -- (101.0406,135.0602) --
    (89.7517,140.7046);
  \path[color=black,draw=black,line join=miter,line cap=butt,miter limit=4.00,line
    width=0.623pt] (149.5453,99.3256) -- (146.4363,91.9051) -- (143.3273,99.3256);
  \path[color=black,draw=black,line join=miter,line cap=butt,miter limit=4.00,line
    width=0.623pt] (53.5851,101.4966) -- (50.4761,94.0761) -- (47.3671,101.4966);

\end{tikzpicture}

{\footnotesize
\caption{Here we draw the construction of the sequence $S_1$ for the example in Figure \ref{orbita1}.}} \label{orbita2} 
\end{figure}

Construct $S_1 = (s_{1,0},\ldots,s_{1,m})$ as follows: $s_{1,0}=r_1$, and $s_{1,j}=\rho(w_j)(r_1)$ as long as $s_{1,j-1}\in (-\infty,d_i)$ if $a_j=x_i$ or $s_{1,j-1}\in (-\infty,f_i(d_i)]$ if $a_j=x_i^{-1}$. We get to $s_{1,k}$, the last element we can define by that process. We must have $k<m-n$: otherwise $d_w(S(\rho,w,r_1),i)\leq d_i$ for some $i$, which is not possible since $d_w(S(\rho,w,x),i)$ is increasing on $x$ (since it is a maximum of increasing homeomorphisms). Choose $s_{1,k+1}> max\{S(\rho,w,q)\}$, and then set $s_{1,j+1}=s_{1,j}+1$ for every $j\geq k+1$. (See figure \ref{orbita2}).

Notice the sequence $S_1$ defines maps $g_i$ on the sets $D_w(S_1,i)$, by taking $g_i(s_{1,j-1})=s_{1,j}$ if $a_j=x_i$ and $g_i(s_{1,j})=s_{1,j-1}$ if $a_j=x_i^{-1}$. Define each $g_i$ on $(-\infty,d_i]\cup D_w(S_1,i)$ so that it agrees with $f_i$ on $(-\infty,d_i]$.

\vs 

{\bf Claim:} The maps $g_i:(-\infty,d_i]\cup D_w(S_1,i)\to\R$ are increasing.

\vs
{\bf Proof:} It is clear that it is increasing on $(-\infty,d_i]$. Notice next that $D_w(S_1,i)\setminus(-\infty,d_i] \subseteq \{s_{1,k},...,s_{1,m}\}$. This makes easy to check that $g_i$ is increasing on this set.

It only remains to show that if $s_{1,l}\in D_w(S_1,i)$ with $l\geq k$ then $g_i(s_{1,l})>g_i(d_i)$. We will distinguish two cases:
\begin{itemize}
\item Case A: $a_{l+1}=x_i$.

By construction of $S_1$ we have $g_i(s_{1,l})=s_{1,l+1}\geq s_{1,k+1}>max\{S(\rho,w,q)\}\geq f_i(d_i)=g_i(d_i)$ as desired.
\item Case B: $a_{l}=x_i^{-1}$.

If $l>k+1$, then $g_i(s_{1,l})=s_{1,l-1}\geq s_{1,k+1}>max\{S(\rho,w,q)\}\geq f_i(d_i)=g_i(d_i)$.
If $l=k+1$, we notice that $s_{1,k}\notin (-\infty,f_i(d_i)]$: Otherwise, following our construction we would have to set $s_{1,k+1}$ as $f_i^{-1}(s_{1,k})$, but that contradicts the definition of $k$. Therefore $g_i(s_{1,k+1})=s_{1,k}>f_i(d_i)=g_i(d_i)$.

Finally, notice that in this case $l\neq k$: If we suppose that $s_{1,k}\notin (-\infty,d_i]$, then $f_{i}(s_{1,k})=s_{1,k-1}\notin (-\infty,f_i(d_i))$ which also contradict the definition of $k$. 
\end{itemize}

This concludes the proof of the claim.
$\hfill \Diamond$
\vs


Next we continue extending the $g_i$.

Now take $r_2$ with $s_{1,m-1}<r_2<s_{1,m}$ and define $S_2 = (s_{2,0},\ldots,s_{2,m})$ by $s_{2,0}=r_2$ and $s_{2,j+1}=s_{2,j}+1$ for $j<m$. As in the previous case, this extends $g_i$ on $D_w(S_2,i)$. To check this extension is increasing observe that $w$ is cyclically reduced, so $a_1\neq a_m^{-1}$. This ensures that there is no problem at the first step, taking $s_{2,0}$ to $s_{2,1}$.

Inductively, construct $S_l = (s_{l,0},\ldots,s_{l,m}) $ from $S_{l-1}$ as we did for $S_2$ from $S_1$. This defines each $g_i$ on $(-\infty,d_i]\cup \bigcup_{l>0} D_w(S_l,i)$, as an increasing and proper map, that agrees with $f_i$ on $(-\infty,d_i]$. Thus each $g_i$ can be extended to $\R$ as an homeomorphism. These extensions can be chosen so that two different $g_{i_1}$ and $g_{i_2}$ do not have any common fixed points after $min\{d_{i_1},d_{i_2}\}$. Therefore $\bigcap_i Fix(g_i) \subset \bigcap_i Fix(f_i)\cap (-\infty,d]$ for $d=min\{d_1,\ldots,d_n\}$.

Since $\overline{\rho}(w)(s_{l,0})=s_{l,m}>s_{l,0}$ we get that $\overline{\rho}(w)$ has no fixed points in $[s_{l,0},s_{l,m}]$ for every $l>0$. Recall that $\overline{\rho}(w)(q)=h(q)>q$, therefore  $\overline{\rho}(w)$ has no fixed points in $[q,+\infty)=[q,h(q)]\cup \ \bigcup_{l>0} [s_{l,0},s_{l,m}]$.

Notice that $\overline{\rho}(w)$ and $h$ are weakly-conjugated by the identity, which is strong on $(-\infty,q]$.

\vsp

{\bf Case Ib:}  $h(q)<q$ and $Fix(h)\cap(q,+\infty)=\emptyset$.

Note that $S(\rho,w,q)=S(\rho,w^{-1},h(q))$. Therefore, we can exchange $w$, $h$ and $q$ by $w^{-1}$, $h^{-1}$ and $h(q)$, and repeat the construction in Case Ia.

\vsp
{\bf  Case IIa:} $h(q)>q$ and $Fix(h)\cap(q,+\infty)$ consists of a single point where $h$ transverses the diagonal.

Repeat the process in case Ia to construct $\rho'\in Rep(\F_n,Homeo_{+}(\R))$ so that for each $i$, $\rho'(x_i)$ and $\rho(x_i)$ agree on $(-\infty,d_i]$ and with $\rho'(w)(x)>x$ for every $x>q$. Moreover, we can ask each $\rho'(x_i)$ to be piecewise linear on $(d_i,+\infty)$. Next we will define $\overline{\rho}$ as a perturbation of $\rho'$ so that $\overline{\rho}(w)$ is conjugated to $h$.

  Take $z>q$ so that every point in $S(\rho',w,z)$ is bigger than $max\{d_i\}$. Let $d_i'=d_{w}(S(\rho',w,z),i)$. We take $z'>max\{d_i'\}+m$ and define $S'=(z',z'-1,...,z'-m)$. Define $\rho''$ so that $\rho''(x_i)$ is a piecewise linear interpolation that extends $\rho'(x_i)_{|(-\infty,d_i']}$ and the map defined on $D_w(S',i)$ as in Case Ia.

  Let $\overline{h}=\rho''(w)$. Since $\overline{h}(z)>z$ and $\overline{h}(z')<z'$ we see that $\overline{h}$ must have a fixed point in $(z,z')$. Let $y=min\{Fix(\overline{h})\cap (z,z')\}$. For our argument we will need $y$ to be a transverse fixed point of $\overline{h}$. If it is not, we will perform an additional perturbation that we turn to describe now.

  Assume $y$ is not transverse. Then $y$ must be a break point of $\overline{h}$, and since it is the first fixed point in $(z,z')$ its left derivative is smaller than $1$. Let $S''=S(\rho'',w,y)=(y_0,y_1,...,y_m)$ and $\epsilon>0$ so that the intervals $(y_j-\epsilon,y_j+\epsilon)$ are either disjoint or identical. For $y_j$ in $D_w(S'',i)$ we make a small perturbation of $\rho''(x_i)$ supported on $[y_j,y_j+\epsilon)$ so that $y_j$ is no longer a break point (if it was one). This perturbation remains piecewise linear, introducing two new break points on $(y_j,y_j+\epsilon)$. The new $\overline{h}$ has a transverse fixed point at $y$, since now $y$ is not a break point and $\overline{h}'(y)<1$ as the left derivative at $y$ has not changed.

We must ensure that there are no new fixed points on $(z,y)$. If $\epsilon$ is small enough, $\overline{h}$ does not change on $(y-\delta,y]$ for some $\delta>0$. This is because for $t\in(y-\delta,y]$ we have that $S(\rho'',w,t)$ is disjoint from the supports of the perturbations. On the other hand, the perturbations are $C^0$ and $\overline{h}$ is away from the diagonal on $[z,y-\delta]$, thus they can be made small enough not to introduce new fixed points on $[z,y-\delta]$.

Having produced the representation $\rho''$ so that $\overline{h}=\rho''(w)$ has a transverse fixed point at $y$, we take $y'>y$ so that $y$ is the only fixed point of $\overline{h}$ on $(z,y')$. Next we will proceed as in Case I, redefining each $\rho''(x_i)$ so that $\rho''(w)$ is unchanged on $(-\infty,y']$ and $\rho''(w)(x)<x$ for every $x>y'$. More precisely, we are in the situation of Case Ib and each $\rho''(x_i)$ gets redefined from the point $d_w(S(\rho'',w^{-1},\overline{h}(y')),i)$.

This new $\rho''$ works as $\overline{\rho}$ in the statement of the claim, as we will check now. It is clear from the construction that $h$ and $\overline{\rho}(w)$ are weakly-conjugated by a map that coincides with $Id$ on $(-\infty,q]$ and each $\overline{\rho}(x_i)$ coincides with $\rho(x_i)$ on $(-\infty,d_i]$. We check that $Fix(\overline{\rho})\subseteq Fix(\rho)$: By Case Ia we have $Fix(\rho')\subseteq Fix(\rho)$. The piecewise linear interpolation for $\rho''$ can be performed without introducing any new global fixed points. So can its further perturbation, following Case Ib. 

\vsp

{\bf Case IIb:} $h(q)<q$ and $Fix(h)\cap(q,+\infty)$ consists of a single point where $h$ transverses the diagonal. 

This case is analogous to IIa, also exchanging $w$, $h$ and $q$ by $w^{-1}$, $h^{-1}$ and $h(q)$. $\hfill\Box$

\vs

Let $\bar \rho$ be the representation given by Lemma \ref{lem interno}. By lemma \ref{promotion} there exists $\varphi\in Homeo_{+}(\R)$ such that $\varphi^{-1}\circ h\circ\varphi=\overline{\rho}(w)$ and $\varphi$ equals $Id$ on $(-\infty,q]$. Define $\rho^*$ as $\varphi\circ\overline{\rho}\circ\varphi^{-1}$. Now $\rho^*(w)=h$ and by construction of $q_0$ we obtain that $\rho^*(x_i)$ and $\rho(x_i)$ coincide over $(-\infty,p)$ as desired.

\section{On commutator varieties}


\label{sec proof commutator}


In this section we prove Lemmas \ref{lem commutator semicont} and \ref{lem commutator} from \S \ref{sec reduction}. Our arguments are based on the analysis of the commutator variety
$$\mathcal V_h:=\{(f,g)\in Homeo_+(\R)\times Homeo_+(\R) \mid [f,g]=h\},$$
of a given homeomorphism of the line $h$. Though very simple, the key observation (and main difference with the strategy for proving Lemma \ref{lem amalgama}) is that the equation $[f,g]=h$ can be rewritten as the equation $$fgf^{-1}=hg.$$ This rewriting provides us the insight that ``$f$ is conjugating $g$ to $hg$''. The idea will be to modify $g$ outside a large compact set, in a way that keeps $g$ and $hg$ conjugated by an element close to $f$.

\vs

In order to control conjugacy class of $hg$ we observe that (see \S \ref{sec conjugacy classes} for definitions)
\begin{itemize}
\item $Inc(hg)=\{x\in\R:g(x)>h^{-1}(x)\}$
\item $Decr(hg)=\{x\in\R:g(x)<h^{-1}(x)\}$
\item $Fix(hg)=\{x\in\R:g(x)=h^{-1}(x)\}$
\end{itemize}
Because of this, in the same way as the conjugacy class of $g$ is determined by the combinatorics of its graph's crossings against the diagonal, we think the conjugacy class of $hg$ as the combinatorics of the crossings of the graph of $g$ against the graph of $h^{-1}$.

It will be handy to have

\begin{defi} For $\phi_1$, $\phi_2$ (partial) homeomorphisms of the line, we define the combinatorics of $(\phi_1,\phi_2)$ as $\C(\phi_1,\phi_2)(x):= sign(\phi_1(x)-\phi_2(x))\in  \{1,-1,0\} $.

\vs

\end{defi}
With this language, $Inc(\phi)=\C(\phi,id)^{-1}(1)$, and $\psi$ is a weak-conjugation from $\phi_1$ to $\phi_2$ if and only if $\C(\phi_1,id)=\C(\phi_2,id)\circ\psi$. Observe that if $\psi\in Homeo_{+}(\R)$ then $\C(\psi\phi_1,\psi\phi_2)=\C(\phi_1,\phi_2)$. This implies that $f$ is a weak conjugation from $g$ to $hg$ if and only if $\C(g,id)=\C(g,h^{-1})\circ f$.

\subsection{Proof of Lemma \ref{lem commutator}}
\label{sec lem 3.3}

Let $K=(-\infty, k]$ be a closed proper interval of the line, and suppose $(f,g)$ is a pair satisfying condition $(*_K)$. We will denote $p_{(f,g,K)}$ from condition $(*_K)$ simply by $p$.

We begin by proving the lemma in a simple case, that will play an important role in the general proof.
\vs

{\bf Toy case:} Assume that $f(k)>k$ and that $Fix(g)\cap(k,f(k)]=\emptyset$.
\vs

In this case, the perturbation of $g$ will be supported on $(f(k),+\infty)$ and the perturbation of $f$ on $(k,+\infty)$. Let's assume that $g(f(k))<f(k)$, the complementary case ($g(f(k))>f(k)$) can be treated identically.

Let's focus first on the construction of $g_2$: From the fact that $f$ conjugates $g$ to $hg$ and that $g(x)<x$ for every $x\in(k,f(k)]$, we get that $hg(f(k))\leq f(k)$. This implies that $g(f(k))\leq min\{f(k),h^{-1}(f(k))\}$, and therefore we can define $g_2$ satisfying $g_2(x)<min\{id,h^{-1}\}(x)$ for every $x>f(k)$.

We are ready to build $f_2$: By the construction of $g_2$ there exists $\psi$ a perturbation of $f$ supported on $(k,+\infty)$ that weakly conjugates $g_2$ and $hg_2$. Since the perturbations are supported outside $(-\infty,k)$, we know that $\psi$ is strong on $(-\infty,min\{k,f^{-1}(k)\})$. Therefore we can apply Lemma \ref{promotion} to ``promote" $\psi$ to a conjugation $f_2$ between $g_2$ and $hg_2$ such that ${f_2}_{|(-\infty,k]}=f_{|(-\infty,k]}$.

Let's focus now on the construction of $g_1$: Since $h^{-1}(p)\neq p$ we can pick by continuity a point $p_1>p$ so that $\C(h^{-1},id)_{|[p,p_1]}$ is constant. ($p$ and $p_1$ in the same ``bump" of the graph of $h^{-1}$). Recall that $g(f(k))\leq min\{f(k),h^{-1}(f(k))\}$. So we can define $g_1$ on $[f(k),p)$ so that $g_1(x)< min\{x,h^{-1}(x)\}$ and $g_1(p)=min\{p,h^{-1}(p)\}$. On $(p,p_1)$ we define $g_1$ so that $min\{x,h^{-1}(x)\}<g_1(x)<max\{x,h^{-1}(x)\}$ and $g_1(p_1)=max\{p_1,h^{-1}(p_1)\}$. Finally we can define $g_1$ on $(p_1,+\infty)$ so that $g_1(x)>max\{x,h^{-1}(x)\}$.

Now we build $f_1$: Observe that both $\C(g_1,id)_{|(k,+\infty)}$ and $\C(g_1,h^{-1})_{|(f(k),+\infty)}$ have a single sign change, that is form $-1$ to $+1$. This implies that we can construct $\psi$ a perturbation of $f$ supported on $(k,+\infty)$ such that $\psi$ weakly-conjugates $g_1$ and $hg_1$. Again, applying Lemma \ref{promotion} we finish the construction.

Finally, notice that $Fix(\langle f_i,g_i\rangle)\subseteq Fix(\langle f,g\rangle)$. Indeed, $g_2$ has no fixed points on the support of the perturbation, while $g_1$ has a single fixed point (either $p$ or $p_1$) and that point is not fixed by $h=[f_2,g_2]$. $\hfill\Diamond$

\vsp

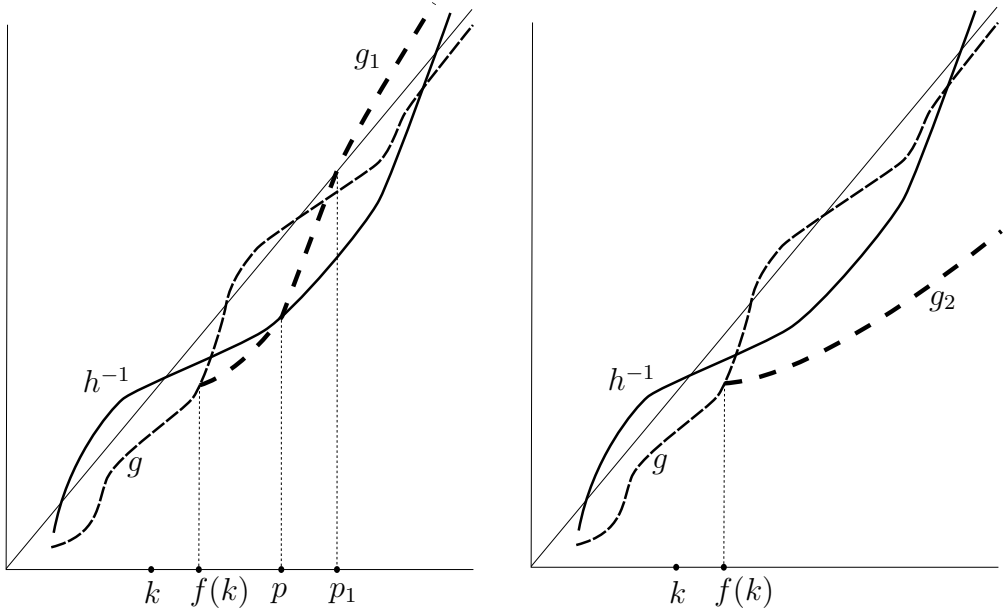
\begin{figure}[h] \label{fig toy case}

\hspace*{-2.7cm}
\begin{tikzpicture}[y=0.80pt, x=0.65pt, yscale=-1.000000, xscale=1.000000, inner sep=0pt, outer sep=0pt]


\begin{scope}[shift={(0,-308.26759)}]
  \path[draw=black,line join=miter,line cap=butt,line width=0.320pt]
    (272.1331,1033.1916) -- (0.2000,1033.1916) -- (1.0081,775.4006);
  \path[draw=black,line join=miter,line cap=butt,line width=0.320pt]
    (0.2013,1032.7888) -- (271.3249,768.1275);
  \path[cm={{0.4,0.0,0.0,0.4,(-15.55838,708.00915)}},fill=black]
    (250.4261,812.9561) ellipse (0.0998cm and 0.1140cm);
  \begin{scope}[shift={(27.71429,0)},shift={(0,0)}]
    \path[cm={{0.4,0.0,0.0,0.4,(-15.55838,708.00915)}},fill=black]
      (250.4261,812.9561) ellipse (0.0998cm and 0.1140cm);
  \end{scope}
  \path[cm={{0.4,0.0,0.0,0.4,(-15.55838,708.00915)}},fill=black]
    (439.7118,812.9561) ellipse (0.0998cm and 0.1140cm);
  \begin{scope}[shift={(32.40296,-0.0002)},shift={(0,0)}]
    \path[cm={{0.4,0.0,0.0,0.4,(-15.55838,708.00915)}},fill=black]
      (439.7118,812.9561) ellipse (0.0998cm and 0.1140cm);
  \end{scope}
  \path[draw=black,dash pattern=on 5.76pt off 0.96pt,line join=miter,line
    cap=butt,miter limit=4.00,line width=0.960pt] (26.4416,1022.6683) .. controls
    (54.1559,1016.9540) and (53.2988,1002.3826) .. (58.4416,990.9540) .. controls
    (63.5845,979.5255) and (103.0130,958.6683) .. (109.0130,951.8112) .. controls
    (115.0130,944.9540) and (125.6900,919.1738) .. (128.1819,907.8587) .. controls
    (130.6737,896.5436) and (139.8702,887.2398) .. (145.8702,881.5255) .. controls
    (151.8702,875.8112) and (207.2988,847.2398) .. (215.5845,841.5255) .. controls
    (223.8702,835.8112) and (227.2988,825.2398) .. (231.2988,818.6683) .. controls
    (235.2988,812.0969) and (271.8702,775.5255) .. (271.8702,775.5255);
  \path[draw=black,dash pattern=on 0.96pt off 0.96pt,line join=miter,line
    cap=butt,miter limit=4.00,line width=0.320pt] (112.3263,1031.5753) --
    (112.7377,946.6643);
  \path[draw=black,line join=miter,line cap=butt,miter limit=4.00,line
    width=0.960pt] (27.8702,1015.8112) .. controls (35.5845,981.8112) and
    (60.1559,958.0969) .. (67.8702,952.6683) .. controls (75.5845,947.2397) and
    (132.7273,929.2397) .. (150.4416,920.0969) .. controls (168.1559,910.9540) and
    (209.0130,871.2398) .. (217.2988,858.3826) .. controls (225.5845,845.5255) and
    (258.7273,769.8112) .. (258.7273,769.8112);
  \path[draw=black,dash pattern=on 0.96pt off 0.96pt,line join=miter,line
    cap=butt,miter limit=4.00,line width=0.320pt] (160.3263,1031.5753) --
    (160.3953,913.6680);
  \path[draw=black,dash pattern=on 0.96pt off 0.96pt,line join=miter,line
    cap=butt,miter limit=4.00,line width=0.320pt] (192.7273,1033.1916) --
    (192.9548,844.6296);
  \path[draw=black,dash pattern=on 6.66pt off 6.66pt,line join=miter,line
    cap=butt,miter limit=4.00,line width=1.664pt] (112.9212,945.9594) .. controls
    (138.7273,939.5255) and (160.8151,913.3582) .. (160.8151,913.3582) .. controls
    (160.8151,913.3582) and (183.7541,863.1348) .. (192.3836,845.1873) .. controls
    (201.0130,827.2397) and (249.2988,765.2397) .. (249.2988,765.2397);
  \path[fill=black] (-121.3352,1188.2791) node[above right] (text4346) {};
  \path[fill=black] (80.2824,1049.2092) node[above right] (text4350) {$k$};
  \path[fill=black] (107.0763,1051.1177) node[above right] (text4358) {$f(k)$};
  \path[fill=black] (154.7101,1049.1022) node[above right] (text4362) {$p$};
  \path[fill=black] (188.4549,1049.1149) node[above right] (text4366) {$p_1$};
  \path[fill=black] (44.9015,949.3335) node[above right] (text4370) {$h^{-1}$};
  \path[draw=black,fill=black,line join=miter,line cap=butt,line width=0.320pt]
    (202.7927,796.6335) node[above right] (text4382) {$g_1$};
  \path[fill=black] (-159.4961,1090.8149) node[above right] (text3784) {};
  \path[fill=black] (70.7049,989.2532) node[above right] (text3788) {$g$};
  \path[draw=black,line join=miter,line cap=butt,line width=0.320pt]
    (577.6597,1032.1120) -- (305.7266,1032.1120) -- (306.5347,774.3210);
  \path[draw=black,line join=miter,line cap=butt,line width=0.320pt]
    (305.7279,1031.7092) -- (576.8515,767.0479);
  \path[cm={{0.4,0.0,0.0,0.4,(289.96824,706.92955)}},fill=black]
    (250.4261,812.9561) ellipse (0.0998cm and 0.1140cm);
  \begin{scope}[shift={(333.24092,-1.0796)},shift={(0,0)}]
    \path[cm={{0.4,0.0,0.0,0.4,(-15.55838,708.00915)}},fill=black]
      (250.4261,812.9561) ellipse (0.0998cm and 0.1140cm);
  \end{scope}
  \path[draw=black,dash pattern=on 5.76pt off 0.96pt,line join=miter,line
    cap=butt,miter limit=4.00,line width=0.960pt] (331.9682,1021.5887) .. controls
    (359.6825,1015.8744) and (358.8254,1001.3030) .. (363.9682,989.8744) ..
    controls (369.1111,978.4459) and (408.5397,957.5887) .. (414.5397,950.7316) ..
    controls (420.5397,943.8744) and (431.2167,918.0942) .. (433.7085,906.7791) ..
    controls (436.2003,895.4640) and (445.3968,886.1602) .. (451.3968,880.4459) ..
    controls (457.3968,874.7316) and (512.8254,846.1602) .. (521.1111,840.4459) ..
    controls (529.3968,834.7316) and (532.8254,824.1602) .. (536.8254,817.5887) ..
    controls (540.8254,811.0173) and (577.3968,774.4459) .. (577.3968,774.4459);
  \path[draw=black,dash pattern=on 0.96pt off 0.96pt,line join=miter,line
    cap=butt,miter limit=4.00,line width=0.320pt] (417.8530,1030.4957) --
    (418.2643,945.5847);
  \path[draw=black,line join=miter,line cap=butt,miter limit=4.00,line
    width=0.960pt] (333.3968,1014.7316) .. controls (341.1111,980.7316) and
    (365.6825,957.0173) .. (373.3968,951.5887) .. controls (381.1111,946.1601) and
    (438.2540,928.1601) .. (455.9683,919.0173) .. controls (473.6825,909.8744) and
    (514.5396,870.1602) .. (522.8254,857.3030) .. controls (531.1111,844.4459) and
    (564.2539,768.7316) .. (564.2539,768.7316);
  \path[draw=black,dash pattern=on 6.66pt off 6.66pt,line join=miter,line
    cap=butt,miter limit=4.00,line width=1.664pt] (418.3237,945.0979) .. controls
    (478.8254,941.8744) and (579.3968,873.0173) .. (579.3968,873.0173);
  \path[fill=black] (385.9849,1049.1465) node[above right] (text3878) {$k$};
  \path[fill=black] (412.2512,1051.4067) node[above right] (text3882) {$f(k)$};
  \path[fill=black] (350.4281,948.2539) node[above right] (text3894) {$h^{-1}$};
  \path[draw=black,fill=black,line join=miter,line cap=butt,line width=0.320pt]
    (537.3968,910.7316) node[above right] (text3898) {$g_2$};
  \path[fill=black] (376.2315,988.1736) node[above right] (text3906) {$g$};
\end{scope}

\end{tikzpicture}
\vspace*{-4cm}
\caption{Toy Case}
\end{figure}

Observe that the Toy Case is analogous to the case where $f(k)<k$ and $Fix(hg)\cap(f(k),k]=\emptyset$. The case $f(k)=k$ is even simpler.

In general, since $f$ is a conjugation from $g$ to $hg$, we know that the combinatorial information of $g$ on $(-\infty,k]$ coincides with the combinatorial information of $hg$ on $(-\infty,f(k)]$, that is
$$\C(g,id) (x) =\C(g,h^{-1})\circ f (x) \text{ for  } x\leq k.$$

The Toy Case was easy because we assumed that the combinatorics of $g$ on $[k,f(k)]$ was constant. In general, we will need to make a previous perturbation in order to attain a similar situation.

\begin{claim} \label{local perturbation} (Local perturbation) {\em There is  $q>p$,  a homeomorphism $\psi:(-\infty,q]\to (-\infty,q]$, and a homeomorphism over its image $\bar g: (-\infty,q]\to \R$ such that
\begin{enumerate}
\item $\C(\bar g,id)=\C(\bar g,h^{-1})\circ\psi$, that is, $\psi $ is a weak conjugation from $\bar{g}$ to $h\bar{g}$ on $(\-\infty,q]$.
\item the pair $(\psi,\bar g)$ agrees with $(f,g)$ on $(-\infty,k]$.
\item the set of common fixed points of $\psi$ and $\bar g$ on $(-\infty,q]$ is contained in the set of common fixed points of $f$ and $g$ on $(-\infty,q]$.

\end{enumerate}
}
\end{claim}

{\bf Proof of the Claim:} If $f(k)=k$, we set $q=k$ and the proof is trivial.
We will distinguish the two cases when $f(k)<k$ and when $f(k)>k$. We will first focus on the construction of the local perturbations $\psi$ and $\bar g$ from the claim, and in the end the third condition of it will be checked.

\vs

{\bf Case I:} $ f(k)<k.$

If $Fix(hg)\cap (f(k),k]=\emptyset$ we proceed as in the Toy Case. Suppose that $Fix(hg)\cap (f(k),k]\neq\emptyset$. Let $p$ be as in the $(*_K)$ condition for the pair $(f,g)$. (For this case we will only need to know that $p>g(k)$, $p>k$ and $h^{-1}(p)\neq p$, that is guaranteed by $(*_K)$).

\vsp

{\bf Subcase Ia:} $ f(k)\notin Fix(hg).$ 

We start by defining $\bar g$ over $(-\infty,p]$. We set $\bar g=g$ on $(-\infty,k]$ and then we extend it over $[k,p]$ satisfying $Fix(\bar g)\cap (k,p]=p$. This is possible because $k\notin Fix(g)$. Let $s_1=min \mbox{ }Fix (hg)\cap (f(k),p]$ and $s_2=max\text{ }Fix(hg)\cap [f(k),p]$. Choose $\epsilon>0$ and define $\psi$ over $(-\infty,p+\epsilon]$ satisfying: $\psi_{|(-\infty,k]}=f_{|(-\infty,k]}$, $\psi(p)=s_1$ and $\psi(p+\epsilon)=s_2$. Now we continue extending $\bar g$. Define $\bar g$ over $[p,p+\epsilon]$ as $\psi^{-1} hg\psi$. Notice that $\bar g$ takes $[p,p+\epsilon]$ to itself. Since $h^{-1}p\neq p$ we can take $\epsilon$ small enough so that $Fix(h\bar g)\cap [p,p+\epsilon]=\emptyset$ (i.e. the graph of $\bar g$ does not meet that of $h^{-1}$ over $[p,p+\epsilon]$). Take $q>p+\epsilon$ and define $\bar{g}$ over $[p+\epsilon,q]$ satisfying $Fix(\bar g)\cap [p+\epsilon,q]=\{p+\epsilon\}$, $Fix(h\bar g)\cap [p+\epsilon,q]=\emptyset$ and $\C(\bar g, id)(q)=\C(\bar g,h^{-
1})(q)$. Finally, any extension of $\psi$ to $(-\infty,q]$ such that $\psi(q)=q$ will satisfy the Claim's thesis.

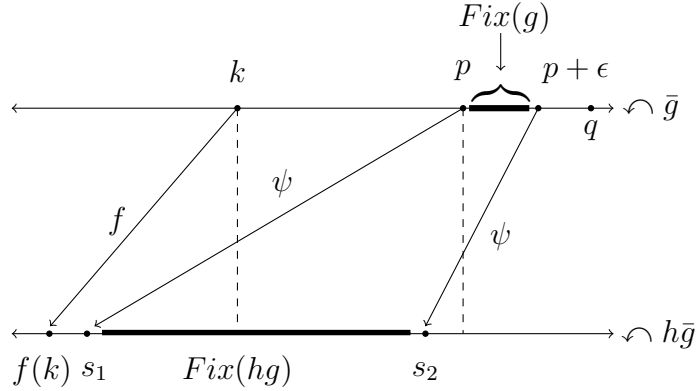
\begin{figure}[h] 
\centering
\begin{tikzpicture}

\draw[<->] (0,3)->(8,3);
\draw (8.5,3) node{$\curvearrowleft \bar{g}$};
\draw[<->] (0,0)->(8,0);
\draw (8.6,0) node{$\curvearrowleft h\bar{g}$};

\filldraw [black] (7.7,3) circle (1 pt);
\draw (7.7,2.7) node{$q$};
\filldraw [black] (7,3) circle (1 pt);
\draw (7.5,3.5) node{$p+\epsilon$};
\filldraw [black] (6,3) circle (1 pt);
\draw (6,3.5) node{$p$};
\filldraw [black] (3,3) circle (1 pt);
\draw (3,3.5) node{$k$};

\draw[very thick,-] (6.08,3.02) -- (6.88,3.02);
\draw[very thick,-] (6.08,2.98) -- (6.88,2.98);
\draw (6.55,4.2) node{$Fix(g)$};
\draw[->]  (6.5,4) -> (6.5,3.5);
\draw (6.5,3.2) node{$\overbrace{}$};

\draw[very thick,-] (1.2,0.03) -- (5.3,0.03);
\draw[very thick,-] (1.2,0) -- (5.3,0);
\draw (3,-0.5) node{$Fix(hg)$};

\draw[dashed] (3,3)--(3,0);
\draw[dashed] (6,3)--(6,0);
\draw[->] (3,3) -> (0.51,0.1);
\draw (1.4,1.5) node{$f$};
\draw[->] (6,3) -> (1.1,0.1);
\draw (3.6,2) node{$\psi$};
\draw[->] (7,3) -> (5.5,0.1);
\draw (6.5,1.3) node{$\psi$};

\filldraw [black] (0.5,0) circle (1 pt);
\draw (0.4,-0.5) node{$f(k)$};
\filldraw [black] (1,0) circle (1 pt);
\draw (1.1,-0.5) node{$s_1$};
\filldraw [black] (5.5,0) circle (1 pt);
\draw (5.5,-0.5) node{$s_2$};

\end{tikzpicture}
\caption{Local perturbation, subcase Ia}
\end{figure}

\vsp

{\bf Subcase  Ib:} $f(k)\in Fix(hg).$ 

Here we begin by defining $\psi$ over $(-\infty,p]$. Take $s\in Fix(hg)\cap (f(k),k]$. Let $\psi:(-\infty,p]\to (-\infty,s]$ be any homeomorphism with $\psi_{|(-\infty,k]}=f_{|(-\infty,k]}$ and $\psi(p)=s$. Now we define $\bar g$ over $(-\infty,p]$ agreeing with $g$ on $(-\infty,k]$ and with $\psi^{-1} hg\psi$ over $(k,p]$. Notice this is well defined since $\psi((k,p])\subseteq (-\infty,k]$. This extension is continuous because $g(k)=k$ and $h\bar g(f(k))=f(k)$. Then we proceed as in subcase Ia, with $s$ as $s_1$.

\vs

{\bf Case II:} $ f(k)>k.$

Since the pair $(f,g)$ satisfies the $(*_K)$ condition, there exists $p>g(f(k))$ such that $h^{-1}(p)\neq p$ and therefore there exists $p_1>p$ such that $h^{-1}(p_1)>p$ and $h^{-1}(p_1)\neq p_1$. ($p_1$ can be taken as an iterate of $p$ by $h^{\pm 1}$). Notice $p_1$ satisfies the conditions for being $p_{(f,g,K)}$ in condition $(*_K)$ and also satisfies $g(f(k))<p<h^{-1}(p_1)$. Thus we can redefine $p:=p_1$ and assume that $g(f(k))<h^{-1}(p)$.

\vsp

{\bf Subcase IIa:} $ k\notin Fix(g).$ 

We begin with the construction of $\bar g$ that will extend $g_{|(-\infty,f(k)]}$. For this define $\bar g$ over $[f(k),p]$ satisfying $Fix(h\bar g)\cap [f(k),p]={p}$. This can be done because $f(k)\notin Fix(h\bar g)$ and $g(f(k))<h^{-1}(p)$.
Consider now $s_1=min\mbox{ }Fix(\bar g)\cap(k,+\infty)$ and $s_2=max\mbox{ }Fix(\bar g)\cap[k,p]$. We now build $\psi$, that will extend $f_{|(-\infty,k]}$. Take $\epsilon>0$ and define $\psi$ on $[k,s_2]$ satisfying $\psi(s_1)=p$ and $\psi(s_2)=p+\epsilon$. We define an auxiliary function $\phi\colon [p,p+\epsilon]\to[p,p+\epsilon]$ as $\phi(x) = \psi\bar g\psi^{-1}(x)$. This is well defined since $\psi^{-1}$ takes $[p,p+\epsilon]$ into $(-\infty,p]$. Then define $\bar g$ over $[p,p+\epsilon]$ as $h^{-1}\phi$. Notice that $\bar g$ is well defined on $p$ since $h^{-1}\phi(p)=h^{-1}(p)$, and we had from before that $p\in\Fix(h\bar g)$. Since $h^{-1}(p)\neq p$, we can show as in case I that if $\epsilon$ is small enough then $Fix(\bar g)\cap(p,p+\epsilon]=\emptyset$. Finally, extend $\bar g$ over $[p+\epsilon,q]$ satisfying $Fix(\bar g)\cap(p+\epsilon,q]=Fix(h\bar g)\cap(p+\epsilon,q]=\emptyset$ and $\C(\bar g,id)(q)=\C(\bar g,h^{-1})(q)$. This allows us to extend $\psi$ over $[s_2,q]$ satisfying the equation $\C(\bar g,id)(x)=\C(\bar g,h^{-1})(\psi(x))$ and with $\psi(q)=q$.

\vsp

{\bf Subcase  IIb:}  $k\in Fix(g).$ 

We can make the construction as in case Ib, with the same modifications we did for case IIa. (Namely, interchanging the roles of $\bar g$ and $h\bar g$).

\vs

Now it only remains to check that $Fix(\langle \psi,\bar g \rangle)\subseteq Fix(\langle f,g \rangle)$. If $x\leq min \ \{k,f(k)\}$ the inclusion is trivial and if $x\in [min \ \{k,f(k)\},q]$ observe that $\psi$ does not fix any point in $Fix(\bar g)$. $\hfill\square$

\vs

Now we show how Lemma \ref{lem commutator} follows from the local perturbation Claim:

Let's begin with the construction of the $g_i$: Suppose that $\C(\bar g,id)(q)=1=\C(\bar g,h^{-1})(q)$. In this case we just define $g_1$ as an extension of $\bar g$ such that $g_1(x)>max\{x,h^{-1}(x)\}$  for every $x\geq q$. To construct $g_2$ we make a small modification on the construction of $\bar g$ and $\psi$ in the Claim's proof. We stop the construction of $\bar g$ at $p+\epsilon$ and then add two transversal intersection points of $\bar g$ with $id$ and with $h^{-1}$, as we did in the Toy Case. (The fact that $h^{-1}(p+\epsilon)\neq p+\epsilon$ implies that we can suppose the two transversal intersections, with $id$ and with $h^{-1}$, occur at different points, so no new global fixed points are created). It's easy to see that we can redefine $\psi$ over $[s_2,q]$ to obtain $\C(\bar g,id)(x)=\C(\bar g,h^{-1})(\psi(x))$ for all $x\leq q$ and $\psi(q)=q$. Now we can ensure that $\C(\bar g,id)(q)=-1=\C(\bar g,h^{-1})(q)$ and define $g_2$ as an extension of $\bar g$ that satisfies $g_2(x)<min \ \{x,h^{-1}(x)\}$ for every $x\geq
q$.

If we had that originally $\C(\bar g,id)(q)=-1=\C(\bar g,h^{-1})(q)$ we proceed analogously.

Now we will construct the $f_i$: observe that we can always extend $\psi$ to the whole line verifying $\C(g_i,id)(x)=\C(g_i,h^{-1})(\psi(x))$ for every $x\in\R$. Therefore we have that $\psi$ is a weak-conjugation between $g_i$ and $hg_i$ that is strong on $(-\infty,min\{k,f^{-1}(k)\}]$. Applying Lemma \ref{promotion} we obtain the $f_i$ extending $f_{|(-\infty,k]}$ and conjugating $g_i$ with $hg_i$ as desired.

It remains to check that $Fix(\langle f_i,g_i \rangle)\subseteq Fix(\langle f,g \rangle)$. If $x\leq q$ it follows from the Local perturbation Claim (even with the modification for $g_2$), because ${f_i}_{|Fix(\bar g)}={\psi}_{|Fix(\bar g)}$. If  $x\geq q$ there are no $g_i$ fixed points. $\hfill\square$

\subsection{Proof of Lemma \ref{lem commutator semicont}}
\label{sec proof of main lemma 3}

We denote by $K'$ the convex closure of $K\cup f(K)$.

First we will assume $K=(-\infty,k]$. As in Lemma \ref{lem commutator}, we distinguish two main cases:

\vs

${\bf Case \mbox{ }I:}\; f(k)<k.$

Here $K'=(-\infty,k]$.

\vsp

$\bf Subcase \mbox{ }Ia:$ {\bf $hg$ has no fixed points in $(f(k),k]$.} Then as in the Toy Case of Lemma \ref{lem commutator} $g(k)$ is not between $k$ and $h^{-1}(k)$. So we can define $g'$ on $(k,+\infty)$ so that it is either less than $min\{id,h'^{-1}\}$ or greater than $max\{id,h'^{-1}\}$ (depending on whether $hg(k)<k$ or $hg(k)>k$). Then $f'$ also follows the construction in Lemma \ref{lem commutator}, and conjugates $g'$ to $h'g'$.

\vsp

$\bf Subcase \mbox{ }Ib:$ {\bf There are fixed points of $hg$ in $(f(k),k]$.} Let $s_1 = min \ Fix(hg)\cap [f(k),k]$ and $s_2 = max \ Fix(hg)\cap [f(k),k]$. If $f(k)<s_1$ we find a point $k_1>k$ so that $k_1>g(k)$, and define $g'$ over $(k,k_1]$ so that $k_1$ is the only fixed point of $g'_{|(k,k_1]}$. If $f(k)=s_1$, let $k_1 = k$. Now let $\psi:(-\infty,k_1+1]\to (-\infty,s_2]$ an homeomorphism agreeing with $f$ on $(-\infty,k]$, and with $\psi(s_1)=k_1$ and $\psi(s_2)=k_1+1$. Extend $g'$ over $(k_1,k_1+1]$ as $\psi^{-1}\circ hg\circ\psi$. (Notice $h'$ agrees with $h$ on $(-\infty,s_2]$).

Now, if $h'g'$ has no fixed points in $(s_2,k_1+1]$, proceed as in Subcase Ia. If there are such fixed points let $s_3 = max \ Fix(h'g')\cap (s_2,k_1+1]$, and we will extend $g'$ and $\psi$ over $(k_1+1,k_1+2]$ as follows. $\psi$ will take $(k_1+1,k_1+2]$ homeomorphically onto $(s_2,s_3]$ and define $g'$ over $(k_1+1,k_1+2]$ as $\psi^{-1}\circ h'g'\circ\psi$. (Notice $g'$ was defined already on $(s_2,s_3]$).

We proceed inductively. The process stops if we fall in Subcase Ia in any instance. Otherwise, notice that $s_{n+3}>k_1+n$ ($s_3>k$, $s_4>k_1+1$ and so on, since $h'g'$ has no fixed points in $(s_{n+2},k_1+n)$). So the $\psi$ obtained is a homeomorphism of the line, that weakly conjugates $g'$ and $h'g'$. We finish the construction by Lemma \ref{promotion}.

We need to check that $Fix(\langle f',g'\rangle)\subseteq Fix(\langle f, g\rangle)$ and that $Fix(\langle f',g'\rangle)\cap(k,+\infty)=\emptyset$. If $x>k$ is a fixed point of $g'$ notice that it is not fixed by $\psi$ in the previous construction. Therefore it is not fixed by $f'$ by Lemma \ref{promotion}. On the other hand $(f',g')$ agrees with $(f,g)$ over $(-\infty,k)$.

\vs

${\bf Case \mbox{ }II: } \;f(k)>k.$

Here $K'=(-\infty,f(k)]$.

Follow the same scheme as in Case I, with the following modifications (similar to those in Case II of the Local Perturbation Claim in \S \ref{sec lem 3.3}): $k_1$ will be taken so that $k_1>f(k)$ and $h'^{-1}(k_1)>g(f(k))$ in the case $f(k)$ is not fixed by $hg$ (and as $f(k)$ otherwise).
On each step of the extension as in Subcase Ib, we extend $h'g'$ over $[k_1+n,k_1+n+1]$ as $\phi=\psi^{-1}\circ g'\circ\psi$ and define $g'=h'^{-1}\phi$ over $[k_1+n,k_1+n+1]$.

It remains the case $f(k)=k$ but it is a simple modification of Subcase Ia.

\marginpar{quedo bien?}
The proof in the case $K=[k,+\infty)$ is analogous. For the case where $K=[u,v]$ is compact, we iterate the case for semi-infinite intervals: Write $K'=[u',v']$ and consider $h_1\in Homeo_{+}(\R)$ that agrees with $h$ on $(-\infty,v')$ and with $h'$ on $[v',+\infty)$. We first apply the lemma for $h_1$ on $(-\infty,v]$. Next we apply it again to the perturbations just obtained for $h'$ on $[u,+\infty)$. $\hfill\square$


\begin{small}


\vspace{1.3cm}

\textit{Juan Alonso}

Fac. Ciencias, Universidad de la Republica Uruguay

juan@cmat.edu.uy

\bigskip

\textit{Joaquin Brum}

Fac. Ingenieria, Universidad de la Republica Uruguay

joaquinbrum@fing.edu.uy

\bigskip

\textit{Crist\'obal Rivas}

Dpto. de Matem\'aticas y C.C., Universidad de Santiago de Chile

cristobal.rivas@usach.cl

\end{small}



\begin{thebibliography}{Dillo 83}




\bibitem{BF} {\sc C. Bonatti, S. Firmo.} Feuilles compact d'un feuilletage générique en codimention 1. {\em Ann. Sci. Éc. Norm. Sup.} {\bf 27} (1994), 407-463.

\bibitem{BG camb} {\sc V. Bludov and A. Glass.}
On free products of of right orderable groups with amalgamated subgroups.
{\em Math. Proc. Camb. Phil. Soc.} {\bf 146} (2009), 591-601.

\bibitem{BG lond} {\sc V. Bludov and A. Glass.}
Word problem, embedding and free products of right-ordered groups with amalgamated subgroups. {\em Proc. Lond. Math. Soc.} {\bf 99} (2009), 585-608.


\bibitem{brodskii} {\sc S. Brodskii.} Equations over groups and groups with one defining relation. {\em Sibirsk Mat. Zh.} {\bf 25} (1984), 84-103. Translation to english in {\em Siberian Math. Journal} {\bf 25} (1984), 235-251.




\bibitem{Calegari} D. Calegari, 
\textit{Circular groups, planar groups, and the Euler class}. 
Geometry and Topology Monographs, Proceedings of the Casson Fest 7 (2004), 431-491.

\bibitem{clay} {\sc A. Clay.}
Free lattice-ordered group and the space of left orderings.
{\em Monatshefte für Mathematik} {\bf 167} (2012), 417-430.

\bibitem{clay rolfsen} {\sc A. Clay, D. Rolfsen.} {\em Ordered groups and topology}.
Graduate Studies in Mathematics 176 (2016).




\bibitem {GOD}{\sc B. Deroin, A. Navas, C. Rivas.} Groups, Orders, and Dynamics. {\em Preprint available on Arxiv} (2015).


\bibitem{GRGA} {\sc B. Farb, D. Fisher (eds).} {\em Geometry, Rigidity, and Group Actions}.
Chicago Lect. Math., Univ. of Chicago Press (2011).


\bibitem {ghys}{\sc É. Ghys.} Groups acting on the circle. {\em Eins. Math.} {\bf 47} (2001), 329-407.


\bibitem {ito alg}{\sc T. Ito.} Dehornoy-like left orderings and isolated left orderings, {\em J. Algebra}, {\bf 374} (2013), 42--58.

\bibitem {ito tohoku}{\sc T. Ito.} Construction of isolated left orderings via partially central cyclic amalgamation,
{\em Tohoku Math. J.} {\bf 68} (2016), 49--71.

\bibitem {ito ggd}{\sc T. Ito.} Isolated orderings on amalgamated free products, {\em Groups. Geom. Dyn}. to appear.


\bibitem{KM} {\sc V. Kopitov, N. Medvedev.} {\em Right-ordered groups}.
Siberian Scholl of Algebra and Logic, Plenum. Publ. Corp., New York (1996).

\bibitem {linnell}{\sc P. Linnell.}  The space of left orders of a groups is either finite or uncountable, {\em Bull. Lond. Math. Soc}. {\bf 43} (2011), 200-202.



\bibitem {mann}{\sc K. Mann.} Spaces of surface group representation. {\em Invent. Math.} {\bf 201} (2015), 669-710.


\bibitem{mann rivas} {\sc K. Mann, C. Rivas.}
Group orderings, dynamics, rigidity.
\textit{Preprint, Available on Arxiv} (2016).

\bibitem{McCoy} {\sc R.McCoy, I. Ntantu.} {\em Topological properties of spaces of continous functions}. Springer Lectures Notes in Math. vol. {\bf 1315}, Springer-Verlag Berlin Heidelberg (1988).



\bibitem{mccleary} S.H. McCleary, 
\textit{Free lattice ordered groups represented as $o$-2 transitive $\ell$-permutation groups}. Trans. Amer. Math. Soc. 209(2) (1985), 69-79.


\bibitem {navas hecke}{\sc A. Navas.} An interesting family of left-ordered groups: Central extensions of Hecke groups. {\em J. Algebra}, {\bf 328} (2011), 32-42.

\bibitem{navas-book} {\sc A. Navas.} {\em Groups of circle diffeomorphisms}.
Chicago Lect. Math., Univ. of Chicago Press (2011).

\bibitem{navas orders} {\sc A. Navas.}
{On the dynamics of left-orderable groups}.
\textit{Ann. Inst. Fourier (Grenoble)} {\bf 60} (2010), 1685-1740.




\bibitem {rivas_free}{\sc C. Rivas.} On the space of left-orderings of free groups. {\em J. Algebra} {\bf 350} (2013), 318-329.


\bibitem {rivas tessera}{\sc C. Rivas, R. Tessera.} On the space of left-orderings of virtually solvable groups. {\em G.G.D.} {\bf 10} (2016),  65-90.


\bibitem{sikora} {\sc A. Sikora.}
Topology on the spaces of orderings of groups.
Bull. Lon. Math. Soc. {\bf 36} (2004): 519-526.

\end{thebibliography}
\end{document}